\documentclass[final,onefignum,onetabnum]{siamart251216}

\usepackage{amsmath, amssymb}
\usepackage{mathrsfs}
\usepackage{bm}
\usepackage{physics}
\usepackage{graphicx}
\usepackage{booktabs}
\usepackage{algorithm}
\usepackage{algpseudocode}
\usepackage{xcolor}
\usepackage{multirow}
\usepackage{mathtools}
\usepackage{tikz}
\usetikzlibrary{shapes.geometric, arrows.meta, positioning, calc, fit, backgrounds, shadows.blur, decorations.pathreplacing}
\usepackage{subcaption}

\definecolor{myblue}{RGB}{60, 130, 200}
\definecolor{mygreen}{RGB}{60, 160, 100}
\definecolor{myred}{RGB}{200, 60, 60}
\definecolor{myorange}{RGB}{220, 140, 40}
\definecolor{mypurple}{RGB}{140, 60, 180}
\definecolor{myteal}{RGB}{40, 160, 160}
\definecolor{bg_stream1}{RGB}{230, 245, 255}
\definecolor{bg_stream2}{RGB}{255, 245, 230}
\definecolor{bg_opt}{RGB}{240, 230, 255}
\definecolor{natureBlue}{RGB}{0, 110, 180}
\definecolor{natureRed}{RGB}{210, 40, 40}
\definecolor{natureGreen}{RGB}{0, 150, 80}
\definecolor{natureGray}{RGB}{120, 120, 120}

\tikzset{
    base/.style = {
        draw,
        rounded corners=3pt,
        align=center,
        font=\small\sffamily,
        blur shadow={shadow blur steps=5}
    },
    startstop/.style = {
        base,
        circle,
        minimum width=1.5cm,
        top color=myred!20,
        bottom color=myred!50,
        draw=myred!80!black,
        font=\bfseries
    },
    process/.style = {
        base,
        rectangle,
        minimum height=1.2cm,
        minimum width=3.5cm,
        top color=white,
        bottom color=gray!20,
        draw=gray!80!black
    },
    int_proc/.style = {
        base,
        rectangle,
        minimum width=4cm,
        top color=myblue!10,
        bottom color=myblue!30,
        draw=myblue!80!black
    },
    mcmc_proc/.style = {
        base,
        rectangle,
        minimum width=4.5cm,
        top color=myorange!10,
        bottom color=myorange!30,
        draw=myorange!80!black
    },
    math_node/.style = {
        base,
        rectangle,
        minimum width=5cm,
        top color=mygreen!10,
        bottom color=mygreen!30,
        draw=mygreen!80!black,
        font=\footnotesize
    },
    opt_node/.style = {
        base,
        rectangle,
        minimum width=5cm,
        top color=mypurple!10,
        bottom color=mypurple!30,
        draw=mypurple!80!black
    },
    decision/.style = {
        base,
        diamond,
        aspect=2,
        minimum width=3cm,
        minimum height=1cm,
        top color=yellow!20,
        bottom color=yellow!50,
        draw=orange!80!black,
        font=\footnotesize\bfseries
    },
    arrow/.style = {thick, -{Stealth[length=3mm]}, rounded corners=5pt},
    data_flow/.style = {thick, -{Stealth[length=3mm]}, dashed, color=myblue, rounded corners=5pt},
    data_flow_stat/.style = {thick, -{Stealth[length=3mm]}, dashed, color=myred, rounded corners=5pt},
    data_flow_z/.style = {thick, -{Stealth[length=3mm]}, dashed, color=myblue, rounded corners=5pt},
    label_text/.style = {font=\scriptsize\bfseries, text=black}
}

\newcommand{\calD}{\mathcal{D}}
\newcommand{\calS}{\mathbb{S}}

\newsiamremark{remark}{Remark}
\newsiamremark{hypothesis}{Hypothesis}
\crefname{hypothesis}{Hypothesis}{Hypotheses}
\newsiamthm{claim}{Claim}
\newsiamremark{fact}{Fact}
\crefname{fact}{Fact}{Facts}
\crefname{appendix}{Appendix}{Appendices}
\Crefname{appendix}{Appendix}{Appendices}

\headers{Projected Sobolev NGD for GPE}{C. Bao, C. Cui, K. Jiang, and S. Shu}

\title{Projected Sobolev Natural Gradient Descent for Efficient Neural Network Solution of the Gross--Pitaevskii Equation\thanks{Submitted to the editors DATE.
\funding{C.~Bao was supported by the National Key R\&D Program of China (2021YFA1001300) and the National Natural Science Foundation of China (12271291).
K.~Jiang was supported by the National Key R\&D Program of China (2023YFA1008802) and the Science and Technology Innovation Program of Hunan Province (2024RC1052).
S.~Shu was supported in part by the National Natural Science Foundation of China (12371373) and the Science Challenge Project (TZ2025007).}}}

\author{Chenglong Bao\thanks{Yau Mathematical Sciences Center, Tsinghua University, Yanqi Lake Beijing Institute of Mathematical
Sciences and Applications, Beijing, China, 100084.}
  \and
  Chen Cui\thanks{Yau Mathematical Sciences Center, Tsinghua University, Beijing, China, 100084.}
  \and
  Kai Jiang\thanks{Hunan Key Laboratory for Computation and Simulation in Science and Engineering,
        Key Laboratory of Intelligent Computing and Information Processing of Ministry
        of Education, School of Mathematics and Computational Science, Xiangtan University, Xiangtan, Hunan,
        China, 411105.}
  \and
  Shi Shu\footnotemark[4]}

\usepackage{amsopn}

\begin{document}

\maketitle

\begin{abstract}
  This paper introduces a projected Sobolev natural gradient descent (NGD) method for computing ground states of the Gross--Pitaevskii equation. By projecting a continuous Riemannian Sobolev gradient flow onto the normalized neural network tangent space, we derive a discrete NGD algorithm that preserves the normalization constraint. The numerical implementation employs variational Monte Carlo with a hybrid sampling strategy to accurately account for the normalization constant. To enhance computational efficiency, a matrix-free Nyström-preconditioned conjugate gradient solver is adopted to approximate the NGD operator without explicit matrix assembly. Numerical experiments demonstrate that the proposed method converges significantly faster than physics-informed neural network approaches and exhibits near-constant wall-clock time and memory over the tested high-dimensional range when the hidden architecture is fixed and the implicit Gram operator is used. Moreover, the resulting neural-network solutions provide high-quality initial guesses that substantially accelerate subsequent refinement by traditional high-precision solvers.
\end{abstract}



\section{Introduction}

Bose--Einstein condensation (BEC) is a macroscopic quantum phenomenon in which many identical bosons occupy a single quantum ground state at ultracold temperatures. Since its first experimental realization in 1995 \cite{anderson1995observation}, BEC has become a central platform for studying macroscopic quantum coherence, superfluidity, and quantum simulation. In the mean-field approximation, the macroscopic wavefunction, or order parameter, is governed by the Gross--Pitaevskii equation (GPE) \cite{lieb2000bosons}. The GPE is a nonlinear Schrödinger equation with a cubic term induced by effective short-range interactions. Computing its ground states is therefore a basic but demanding problem in numerical analysis. Mathematically, the task is to minimize a highly non-convex energy functional on a Riemannian manifold defined by a normalization constraint. The interaction nonlinearity can make the energy landscape severely ill-conditioned, which complicates robust and efficient optimization.

Existing numerical methods for GPE ground states broadly follow two paradigms. The first is discretize--then--solve: spatial discretization reduces the continuous problem to a nonlinear eigenvalue problem, which is then treated by algebraic techniques such as self-consistent field iterations \cite{dion2007ground,cances2021convergence} and inverse iteration schemes \cite{jarlebring2014inverse}. The second is optimize--then--discretize: the continuous energy functional is minimized directly under the manifold constraint. Representative methods include projected Sobolev gradient methods \cite{bao2004computing,garcia2001optimizing,danaila2010new,henning2020sobolev,chen2024convergence}, Riemannian conjugate gradient algorithms \cite{antoine2017efficient,danaila2017computation,peterseim2025energy}, momentum-accelerated schemes \cite{chen2023second}, and preconditioning strategies \cite{antoine2017efficient,feng2025preconditioned}. A comprehensive overview is given in \cite{henning2025gross}.

Despite their maturity, these algorithms can face reduced efficiency and robustness in extreme physical regimes and high-dimensional settings. In the presence of strong interactions or highly heterogeneous external potentials, the GPE energy landscape becomes severely ill-conditioned, which often leads to stagnation or convergence to spurious local minima \cite{peterseim2025neural}. Moreover, although the GPE is conventionally formulated in physical space, the investigation of its high-dimensional extensions is of both physical and numerical relevance. On the physical side, modern research directions such as multicomponent spinor condensates \cite{ueda2022spinor} and quantum Hall physics in synthetic dimensions \cite{ozawa2019topological} naturally give rise to nonlinear Schrödinger equations posed in dimensions exceeding three. From a numerical standpoint, high-dimensional formulations are often necessary for modeling quasi-periodic potentials, where auxiliary dimensions are introduced to restore periodicity \cite{jiang2025irrational}. In these contexts, classical grid-based discretizations are hampered by the curse of dimensionality, as their memory footprint and computational cost increase exponentially with the dimension, rendering them impractical for large-scale simulations.

Deep neural networks (DNNs) have recently become effective tools for scientific computing, especially for high-dimensional approximation problems \cite{raissi2019physics,denis2025accurate,wu2024variational,viteritti2023transformer}. Motivated by these advances, several neural methods for the GPE exploit its structural and geometric properties \cite{bao2025computing,kong2025rotating}. For example, the norm-DNN method \cite{bao2025computing} augments physics-informed neural networks (PINNs) \cite{raissi2019physics} with an explicit normalization layer. However, direct PINN-type formulations for the GPE face two main challenges. First, they typically evaluate residual-based losses on fixed or uniformly sampled collocation points. In high dimensions, and especially for strongly localized probability densities, such sampling can miss the physically relevant regions. The resulting gradients can then be small at most sampled points, causing parameter updates to stagnate. Second, training usually relies on first-order optimizers such as Adam \cite{kingma2015adam}, which can converge slowly on ill-conditioned high-dimensional energy landscapes.

To circumvent the sampling limitations inherent in collocation-based schemes, neural network representations of wavefunctions, commonly termed neural quantum states (NQS), have emerged as a compelling alternative. In the context of linear quantum many-body physics, NQS combined with neural variational Monte Carlo (VMC) methods \cite{carleo2017solving} have achieved state-of-the-art performance, with notable successes including FermiNet \cite{pfau2020ab} and PauliNet \cite{hermann2020deep}. In contrast to PINNs, VMC naturally handles high-dimensional integration via importance sampling from the wavefunction-induced distribution and structurally enforces normalization. However, extending this paradigm to the GPE presents distinct difficulties. Unlike the linear Schrödinger equation, the GPE includes a nonlinear interaction term that complicates the construction of low-variance gradient estimators.

From an optimization perspective, several high-performance neural optimization methods have been developed, including evolutional DNNs \cite{du2021evolutional,kim2025multi,kao2025petnns,zhang2025low}, neural Galerkin methods \cite{bruna2024neural,berman2023randomized,li2025mad}, and natural gradient descent (NGD) methods \cite{martens2020new,zeinhofer2024unified,bonfanti2024challenges,nurbekyan2023efficient,chen2024teng,schwencke2025amstramgram,guzman2025improving,dangel2024kronecker}. These approaches are closely related to the Dirac--Frenkel variational principle (DFVP) \cite{dirac1930note,frenkel1934wave}, which projects continuous-time dynamics onto finite-dimensional parameter manifolds. Existing studies have examined metric selection, temporal discretization, and efficient approximations of the associated metric matrices. Their direct application to the GPE nevertheless remains nontrivial. The GPE combines strong nonlinearity with a strict normalization constraint, requiring a variational optimization scheme tailored to this structure.

In this work, we propose a projected Sobolev NGD method within a neural VMC framework for computing ground states of high-dimensional GPEs. Rather than treating neural-network training as a purely Euclidean optimization problem in parameter space, we use an \emph{optimize--then--discretize} strategy based on the geometry of the normalized function manifold. Starting from a projected Sobolev gradient flow, we derive natural-gradient dynamics in parameter space through the DFVP \cite{dirac1930note,frenkel1934wave}. This construction preserves the geometric constraints and links continuous functional optimization to discrete parameter updates, as illustrated in \Cref{fig:PSNGD}.
To address strong nonlinearity and high dimensionality, we further develop a set of implementation techniques that target the main computational bottlenecks. The hybrid sampling strategy combines Markov chain Monte Carlo (MCMC) with adaptive defensive importance sampling (ADIS) to estimate the energy gradient, normalization constant, and Sobolev Gram quantities. The matrix-free Nyström-preconditioned conjugate gradient solver avoids explicit Gram-matrix assembly and inversion, thereby reducing memory use and improving scalability.
The numerical results show three advantages of the proposed framework. First, compared with existing approaches such as norm-DNN, the method converges faster, often by an order of magnitude (see \Cref{tab:normdnn_comparison}). Second, the optimized neural-network solution provides a high-quality initial guess for traditional high-precision solvers, reducing both iteration counts and time to convergence (see \Cref{fig:initialization}). Finally, the implicit Gram-operator implementation substantially reduces dimension dependence in the high-dimensional experiment: with a fixed hidden architecture, wall-clock time and GPU memory remain nearly constant for $d=4,\ldots,8$, with no evidence of exponential growth (see \Cref{tab:hd_verify}).

The remainder of this paper is organized as follows. \Cref{sec:variational} reviews the variational formulation of the GPE and derives the continuous projected Sobolev gradient flow. \Cref{sec:neural_discretization} introduces the neural discretization and derives the resulting NGD scheme. \Cref{sec:vmc} presents the VMC implementation, including gradient estimation, adaptive importance sampling for the normalization constant, and randomized preconditioning of the Gram matrix. \Cref{sec:numerical_experiments} reports numerical experiments assessing the proposed method. \Cref{sec:conclusion} concludes with limitations and future directions. The appendices collect the supporting stability analysis and ablation studies.

\section{Variational formulation and Sobolev gradient flow}
\label{sec:variational}

This section sets up the variational framework for the GPE and derives the continuous projected Sobolev gradient flow used for the neural discretization.

\subsection{Gross--Pitaevskii energy functional and ground state problem}
\label{subsec:gpe-functional}

The dimensionless GPE on a bounded domain $\calD \subset \mathbb{R}^d$
with Lipschitz boundary $\partial\calD$ describes
a BEC in the mean-field approximation at zero temperature. Its energy functional is
\begin{equation}
  \label{eq:standard_energy_functional}
  E(\psi) = \int_{\calD} \left( \frac{1}{2} |\nabla \psi(\mathbf{r})|^2
    + V(\mathbf{r}) |\psi(\mathbf{r})|^2
  + \frac{\beta}{2} |\psi(\mathbf{r})|^4 \right) \,\mathrm{d}\mathbf{r},
\end{equation}
where $\psi : \calD \to \mathbb{R}$ is the condensate wavefunction subject to the
mass constraint $\|\psi\|_{L^2(\calD)}^2 = 1$. The parameter $\beta \in \mathbb{R}$
represents the dimensionless interaction strength, where $\beta>0$ corresponds to
the repulsive interactions investigated in this work and $\beta<0$ models attractive
interactions. The external potential $V(\mathbf{r})$ is assumed to be
real-valued, essentially bounded, i.e., $V\in L^\infty(\calD)$, and $V(x) \geq 0$ almost everywhere.
We impose homogeneous Dirichlet boundary conditions $\psi=0$ on $\partial\calD$
and use $H_0^1(\calD)$ as the basic Sobolev space associated with this boundary
condition.

The ground state is defined as the global energy minimizer on the unit $L^2$ sphere,
\begin{equation}
  \label{eq:minimization}
  \psi_{\mathrm{GS}}
  = \underset{\psi \in \calS}{\arg\min}\, E(\psi),
  \qquad
  \calS := \left\{ \psi \in H_0^1(\calD) \cap L^4(\calD) :
  \|\psi\|_{L^2(\calD)}^2 = 1 \right\}.
\end{equation}
The intersection with $L^4(\calD)$ defines the finite-energy admissible class for the quartic interaction and ensures that $\int_{\calD}|\psi|^4\,\mathrm{d}\mathbf{r}$ is finite for every $\psi\in\calS$. This condition is automatic for $d\le 4$ through the Sobolev embedding $H_0^1(\calD)\hookrightarrow L^4(\calD)$ and becomes an effective restriction only when $d>4$. Critical points of $E$ on $\calS$ correspond to stationary states and satisfy a constrained variational principle. By introducing a Lagrange multiplier $\mu\in\mathbb{R}$, these states satisfy $\delta(E(\psi) - \mu(\|\psi\|_{L^2}^2 - 1))=0$. The first variation in a direction $v\in H_0^1(\calD)\cap L^4(\calD)$ is given by
\begin{equation}
  \langle E'(\psi), v \rangle
  = \int_{\calD} \left(
    \nabla \psi \cdot \nabla v
    + 2V \psi v
    + 2\beta |\psi|^2 \psi v
  \right)\,\mathrm{d}\mathbf{r},
  \label{eq:fretchet_derivative}
\end{equation}
leading to the nonlinear eigenvalue problem
\begin{equation}
  H[\psi]\psi :=
  \left( -\frac{1}{2}\Delta + V(\mathbf{r}) + \beta |\psi|^2 \right)\psi
  = \mu \psi.
  \label{eq:NSEVP}
\end{equation}
For a ground state $\psi_{\mathrm{GS}}$, the associated chemical potential satisfies $\mu = E(\psi_{\mathrm{GS}}) + \frac{\beta}{2}\|\psi_{\mathrm{GS}}\|_{L^4(\calD)}^4$. The ground state can be chosen nonnegative in $\calD$.

\subsection{Projected Sobolev gradient flow and metric selection}
\label{sec:continuous_sobolev_gradient}

To compute ground states while maintaining the mass constraint, we give a formal Riemannian projected Sobolev-gradient-flow derivation \cite{henning2025gross} that integrates metric preconditioning with projection onto the mass-orthogonal admissible directions. At $\psi$, these directions are
\[
  T_\psi\calS
  =
  \bigl\{
    v\in H_0^1(\mathcal D)\cap L^4(\mathcal D)
    :
    (\psi,v)_{L^2}
    =
    \int_{\mathcal D}\psi v\,\mathrm d\mathbf r
    =0
  \bigr\}.
\]
To state the Riesz construction rigorously, one must first select a Hilbert space $X$ compactly embedded in $L^2(\mathcal D)$, compatible with the admissible test class, and such that $E'(\psi)$ is continuous on $X$. Whenever such a Hilbert structure has been fixed, the Riesz representer of $E'(\psi)$ defines the $X$ \emph{Sobolev gradient} $\nabla_X E(\psi)\in X$ through
\begin{equation}
  (\nabla_X E(\psi),\phi)_X
  =
  \langle E'(\psi),\phi\rangle,
  \qquad
  \forall \phi\in X,
  \label{eq:sobolev_gradient_def_revised}
\end{equation}
where $\langle\cdot,\cdot\rangle$ denotes the duality pairing.

Formally choosing $X=L^2(\mathcal D)$ gives the $L^2$ gradient $\nabla_{L^2} E(\psi) = 2H[\psi]\psi$, and the resulting flow recovers the imaginary time evolution method \cite{bao2004computing}. However, this approach often suffers from numerical stiffness because $H[\psi]$ is an unbounded operator on $L^2$. To mitigate this issue, we use an energy-adaptive metric \cite{henning2020sobolev,henning2025metric} induced by the bilinear form
\begin{equation}
  a_\psi(v,w)
  :=
  \int_{\mathcal D}
  \Bigl(
    \tfrac12\nabla v\cdot\nabla w
    +
    Vvw
    +
    \beta|\psi|^2vw
  \Bigr)\,
  \mathrm d\mathbf r,
  \label{eq:a_inner}
\end{equation}
evaluated on $H_0^1(\mathcal D)\cap L^4(\mathcal D)$. For $d>4$, the norm induced by $a_\psi$ does not in general control the $L^4$ norm, so this intersection need not be complete in the $a_\psi$ norm. We therefore do not assert that $a_\psi$ makes $H_0^1(\mathcal D)\cap L^4(\mathcal D)$ a Hilbert space. A fully rigorous Hilbert-manifold treatment in arbitrary dimension would require an additional choice of topology for which the relevant embedding and the continuity and coercivity properties of $a_\psi$ are established. In the present work, the following continuous-space construction is understood formally on the sufficiently regular finite-energy class represented by the neural ansatz. On this class, we write the formal $a_\psi$-Sobolev gradient as the function satisfying
$a_\psi(\nabla_{a_\psi}E(\psi),\phi) = \langle E'(\psi),\phi\rangle$ for all sufficiently regular finite-energy test functions $\phi$. With the convention in \eqref{eq:fretchet_derivative}, $\langle E'(\psi),\phi\rangle = 2a_\psi(\psi,\phi)$. Thus, the metric incorporates leading-order curvature information from the energy functional and acts as a nonlinear, state-dependent preconditioner for the $L^2$ flow.

Because the normalization constraint is imposed in the $L^2$ geometry, projecting an $X$-gradient onto $T_\psi\mathbb S$ requires the $X$-Riesz representer $n_\psi\in X$ of the constraint differential, defined by $(n_\psi,\phi)_X = (\psi,\phi)_{L^2}$ for all $\phi\in X$. The Riemannian projection is defined as
\begin{equation}
  \mathrm{grad}_X E(\psi)
  =
  P_{\psi,X}(\nabla_X E(\psi))
  =
  \nabla_X E(\psi)
  -
  \frac{(\nabla_X E(\psi),n_\psi)_X}{\|n_\psi\|_X^2}\,n_\psi.
  \label{eq:projection_operator_revised}
\end{equation}
The continuous projected Sobolev gradient flow is then given by
\begin{equation}
  \partial_t \psi(t)
  =
  -\,\mathrm{grad}_X E(\psi(t)),
  \label{eq:projected_gradient_flow_revised}
\end{equation}
which preserves the mass $\|\psi(t)\|_{L^2}^2\equiv1$ and dissipates the energy monotonically.

\begin{remark}
  Beyond formulations based on first-order ordinary differential equations (ODEs), various alternative computational strategies exist for obtaining BEC ground states. These include second-order ODE systems analogous to Nesterov momentum \cite{chen2023second}, Riemannian conjugate gradient methods \cite{danaila2017computation}, and preconditioned normalized $L^2$ flows \cite{antoine2017efficient}. Although these approaches differ in implementation complexity and convergence behavior, each could serve as a starting point for deriving discrete neural schemes similar to the one developed here.
\end{remark}

\section{Neural discretization and natural gradient descent}
\label{sec:neural_discretization}

This section introduces the neural network discretization of the GPE. By invoking the DFVP, we project the continuous Sobolev gradient flow onto the parameterized manifold, which yields a natural gradient evolution for the network parameters. We first detail the neural network architecture and normalization strategies before deriving the discrete parameter-space dynamics.

\subsection{Deep neural network ansatz and normalization}
\label{sec:nn_ansatz}

Since the ground state can be chosen nonnegative, the implementation parameterizes its log amplitude by combining a Gaussian envelope with a multilayer perceptron (MLP) correction. For an input $\mathbf r\in\mathbb{R}^d$, the forward propagation of the $L$-layer correction network is
\begin{align}
  \mathbf{h}^{(0)} &= \mathbf{r}, \nonumber\\
  \mathbf{h}^{(l)} &= \sigma\!\left(\mathbf{W}^{(l)}\mathbf{h}^{(l-1)} + \mathbf{b}^{(l)}\right),
  \qquad l=1,\dots,L-1, \\
  f_{\boldsymbol\eta}(\mathbf{r}) &= \mathbf{W}^{(L)}\mathbf{h}^{(L-1)} + \mathbf{b}^{(L)}, \nonumber
\end{align}
where $\sigma$ denotes the activation function and $\boldsymbol\eta=\{\mathbf{W}^{(l)},\mathbf{b}^{(l)}\}_{l=1}^L$ contains the MLP parameters. The unnormalized wavefunction used in the experiments is
\begin{equation}
  \psi_{\boldsymbol\theta}(\mathbf r)
  =
  \exp\!\left(-\tfrac12\|\mathbf C^\top\mathbf r\|^2\right)
  \exp\!\left(f_{\boldsymbol\eta}(\mathbf r)\right),
  \qquad
  \boldsymbol\theta=(\boldsymbol\eta,\mathbf C),
\end{equation}
where the trainable matrix $\mathbf C$ provides a Gaussian envelope and the MLP learns a multiplicative correction. The experimental implementation registers a full $d\times d$ trainable matrix and uses its lower-triangular part in the forward evaluation.

A key challenge in learning GPE wavefunctions is enforcing the normalizing constraint while optimizing unconstrained network parameters. Standard VMC methods for linear Hamiltonians can work with unnormalized densities by exploiting scale invariance \cite{pfau2020ab}. This invariance is lost for the GPE because of its nonlinear interaction term.

One way to impose the constraint is to use the normalized ansatz \cite{bao2025computing,kong2025rotating}
\begin{equation}
  \hat{\psi}_{\boldsymbol{\theta}}(\mathbf{r}) = \frac{\psi_{\boldsymbol{\theta}}(\mathbf{r})}{\|\psi_{\boldsymbol{\theta}}\|_{L^2(\mathcal{D})}}.
\end{equation}
However, differentiating the normalization integral by automatic differentiation is computationally expensive and can be unstable. We therefore use a hybrid strategy: the normalized ansatz enforces mass conservation, while the log-derivative statistical estimators described in \Cref{subsec:grad} compute gradients without directly differentiating the normalization factor.

\subsection{Natural gradient descent}
\label{sec:natgrad}

We introduce the normalized neural manifold $\mathcal M_{\boldsymbol{\theta}}= \{\hat\psi_{\boldsymbol{\theta}}\in\mathbb S:\boldsymbol{\theta}\in\mathbb R^P\}$. To project the gradient flow \eqref{eq:projected_gradient_flow_revised} onto this manifold, we use the chain rule, which gives
$
\partial_t\hat\psi_{\boldsymbol{\theta}} = \mathbf{J}(\boldsymbol{\theta})\dot{\boldsymbol{\theta}},
$
where $\mathbf{J}(\boldsymbol{\theta})=[\partial_{\theta_1}\hat\psi_{\boldsymbol{\theta}},\dots,\partial_{\theta_P}\hat\psi_{\boldsymbol{\theta}}]$ denotes the Jacobian operator. Because $\hat\psi_{\boldsymbol{\theta}}$ is normalized, its parameter derivatives obey
$
\langle\hat\psi_{\boldsymbol{\theta}},\partial_{\theta_j}\hat\psi_{\boldsymbol{\theta}}\rangle_{L^2}=0,
$
so each $\partial_{\theta_j}\hat\psi_{\boldsymbol{\theta}}$ lies in the tangent space $T_{\hat\psi_{\boldsymbol{\theta}}}\mathbb S$. The parameter-induced tangent subspace is therefore
\[
  \mathcal V_{\boldsymbol{\theta}}
  =
  \mathrm{span}\{\partial_{\theta_j}\hat\psi_{\boldsymbol{\theta}}:\,j=1,\dots,P\}
  \subset T_{\hat\psi_{\boldsymbol{\theta}}}\mathbb S.
\]

\begin{remark}
\label{rem:pseudotime}
Here $t$ is only a pseudo-time indexing the optimization path $\boldsymbol{\theta}(t)$, not the physical time in the time-dependent GPE.
\end{remark}

Since the neural ansatz possesses finite expressivity, the exact Riemannian gradient $\mathrm{grad}_X E(\hat\psi_{\boldsymbol{\theta}})$ generally extends beyond the tangent subspace $\mathcal V_{\boldsymbol{\theta}}$. The DFVP \cite{dirac1930note, frenkel1934wave} resolves this by seeking a velocity $\dot{\boldsymbol{\theta}}$ whose dynamical residual is orthogonal to $\mathcal V_{\boldsymbol{\theta}}$ under the $X$-metric. This approach, also known as the Galerkin projection \cite{bruna2024neural,armegioiu2025functional}, yields
\begin{equation}
  \bigl\langle
  \mathbf{J}(\boldsymbol{\theta})\dot{\boldsymbol{\theta}}
  +
  \mathrm{grad}_X E(\hat\psi_{\boldsymbol{\theta}}),
  \partial_{\theta_j}\hat\psi_{\boldsymbol{\theta}}
  \bigr\rangle_X
  =0,
  \qquad
  j=1,\dots,P.
  \label{eq:galerkin_condition_revised}
\end{equation}

Expanding the above equation leads to a system of ODEs governing the evolution of the parameters
\[
  \sum_{i=1}^P
  \dot{\theta}_i
  \bigl(
    \partial_{\theta_i}\hat\psi_{\boldsymbol{\theta}},
    \partial_{\theta_j}\hat\psi_{\boldsymbol{\theta}}
  \bigr)_X
  =
  -
  \bigl(
    \mathrm{grad}_X E(\hat\psi_{\boldsymbol{\theta}}),
    \partial_{\theta_j}\hat\psi_{\boldsymbol{\theta}}
  \bigr)_X,
  \qquad j=1,\dots,P.
\]
To simplify notation, we introduce the Gram matrix associated with the parameterized manifold $\mathcal M_{\boldsymbol{\theta}}$,
\begin{equation}
  \mathbf{G}_{ij}(\boldsymbol{\theta})
  =
  \bigl(
    \partial_{\theta_i}\hat\psi_{\boldsymbol{\theta}},
    \partial_{\theta_j}\hat\psi_{\boldsymbol{\theta}}
  \bigr)_X,
  \label{eq:gram_matrix_revised}
\end{equation}
which is symmetric and positive semidefinite.
At this stage, we keep the inner product $(\cdot,\cdot)_X$ abstract to emphasize that the Galerkin projection and the resulting parameter dynamics depend on the chosen metric but are not tied to a particular inner product. For any such choice, the Gram matrix $\mathbf{G}(\boldsymbol{\theta})$ is the pullback of the corresponding function-space metric under the normalized neural embedding $\boldsymbol{\theta}\mapsto\hat\psi_{\boldsymbol{\theta}}$. The energy-adaptive choice used in this work is specified after the natural-gradient iteration is derived.

Using the definition of the Sobolev gradient in \eqref{eq:sobolev_gradient_def_revised} and the fact that $\partial_{\theta_j}\hat\psi_{\boldsymbol{\theta}}\in T_{\hat\psi_{\boldsymbol{\theta}}}\mathbb S$, the right-hand side can be written entirely in parameter derivatives. Specifically, we obtain
\begin{align*}
  \bigl( \mathrm{grad}_X E(\hat\psi_{\boldsymbol{\theta}}), \partial_{\theta_j} \hat{\psi}_{\boldsymbol{\theta}} \bigr)_X
  &= \bigl( \nabla_X E(\hat\psi_{\boldsymbol{\theta}}), \partial_{\theta_j} \hat{\psi}_{\boldsymbol{\theta}} \bigr)_X \\
  &= \langle E'(\hat{\psi}_{\boldsymbol{\theta}}), \partial_{\theta_j} \hat{\psi}_{\boldsymbol{\theta}} \rangle
  = \frac{\partial}{\partial \theta_j} E(\hat{\psi}_{\boldsymbol{\theta}})
  = \nabla_{\boldsymbol{\theta}} L(\boldsymbol{\theta})_j,
\end{align*}
where we have defined the parameter-space objective function as $L(\boldsymbol{\theta}):=E(\hat{\psi}_{\boldsymbol{\theta}})$.
Therefore, the evolution of the parameters is governed by the natural gradient flow
\begin{equation}
  \mathbf{G}(\boldsymbol{\theta})\,\dot{\boldsymbol{\theta}}
  =
  -\,\nabla_{\boldsymbol{\theta}} L(\boldsymbol{\theta}),
  \label{eq:natural_gradient_flow}
\end{equation}
which constitutes a finite-dimensional realization of the Sobolev gradient flow on the neural network manifold.

Applying the forward Euler method to \eqref{eq:natural_gradient_flow} leads to the following NGD iteration
\begin{equation}
  \boldsymbol{\theta}^{k+1} = \boldsymbol{\theta}^k - \tau_k \, \mathbf{G}(\boldsymbol{\theta}^k)^{-1} \nabla_{\boldsymbol{\theta}} L(\boldsymbol{\theta}^k),
  \label{eq:nat_grad_descent}
\end{equation}
where $\tau_k>0$ denotes the learning rate.
\Cref{fig:PSNGD} illustrates both the algebraic structure and the geometric interpretation of the proposed method.  Panel~(a) depicts a commutative diagram emphasizing manifold consistency: the discrete NGD update on $\Theta$ represents the optimal finite-dimensional approximation of the continuous Sobolev gradient flow in $\calS$. Panel~(b) geometrically interprets the DFVP-induced Galerkin projection $\mathcal{P}$. Here, the NGD update identifies the descent direction within the tangent space $T_{\hat{\psi}}\mathcal{M}_{\boldsymbol{\theta}}$ that minimizes the angular deviation from the true Sobolev gradient under the $X$-inner product.
\begin{figure}[!htbp]
  \centering
  \resizebox{\textwidth}{!}{
    \begin{tikzpicture}[scale=0.9, transform shape, font=\rmfamily]
      \begin{scope}[xshift=-6.8cm, yshift=-0.5cm]
        \node[anchor=south, font=\bfseries\large] at (3.0, 5.2) {(a) Algebraic Correspondence};

        \node (Theta_k) at (0, 0) [circle, fill=natureBlue!5, draw=natureBlue, thick, inner sep=6pt] {$\boldsymbol{\theta}^k$};
        \node (Theta_target) at (6.0, 0) [circle, fill=natureBlue!5, draw=natureBlue, thick, inner sep=6pt] {$\boldsymbol{\theta}^*$};
        \node (Psi_k) at (0, 3.5) [circle, fill=natureRed!5, draw=natureRed, thick, inner sep=6pt] {$\hat{\psi}_{\boldsymbol{\theta}^k}$};
        \node (Psi_target) at (6.0, 3.5) [circle, draw=natureRed, dashed, thick, inner sep=6pt, fill=white] {$\hat{\psi}_{\text{target}}$};

        \node[natureRed, font=\small, anchor=south] at (0, 4.5) {\textbf{Function Space} $\calS$};
        \node[natureBlue, font=\small, anchor=north] at (0, -1.0) {\textbf{Param Manifold} $\Theta$};

        \draw[->, very thick, natureRed, >={Stealth[length=3mm]}] (Psi_k) -- node[above, font=\small, align=center, yshift=2pt] {
          \textbf{Continuous Sobolev Flow}\\
          $\partial_t \hat{\psi} = - \text{grad}_X E(\hat{\psi})$
        } (Psi_target);

        \draw[->, double, double distance=1.5pt, thick, natureGray, >={Stealth[length=3mm]}] (Theta_k) -- node[left, font=\small, align=right, xshift=-4pt] {
          Neural \\Ansatz
        } (Psi_k);

        \draw[->, thick, natureGreen, >={Stealth[length=3mm]}] (Psi_target) -- node[right, font=\small, align=left, xshift=4pt] {
          \textbf{Galerkin}\\ \textbf{Projection}
        } (Theta_target);

        \draw[->, very thick, natureBlue, >={Stealth[length=3mm]}] (Theta_k) -- node[below, font=\small, align=center, yshift=-4pt] {
          \textbf{Natural Gradient Flow}\\
          $\dot{\boldsymbol{\theta}} = - \mathbf{G}(\boldsymbol{\theta})^{-1} \nabla_{\boldsymbol{\theta}} L(\boldsymbol{\theta})$
        } (Theta_target);

        \node at (3.0, 1.75) [font=\LARGE, natureGray] {$\circlearrowright$};
        \node at (3.0, 1.2) [font=\tiny, natureGray] {Commutative};
      \end{scope}

      \begin{scope}[xshift=2.2cm, yshift=0cm]
        \node[anchor=south, font=\bfseries\large] at (2.5, 4.9) {(b) Tangent Space Projection};

        \shade[left color=natureBlue!5, right color=natureBlue!15, opacity=0.6]
        plot [smooth cycle, tension=0.8] coordinates {(-0.5, 1.0) (2.5, -1.5) (6.4, 0.5) (5.0, 4.5) (0.5, 3.5)};
        \draw[thick, natureBlue!60]
        plot [smooth cycle, tension=0.8] coordinates {(-0.5, 1.0) (2.5, -1.5) (6.4, 0.5) (5.0, 4.5) (0.5, 3.5)};

        \node at (5.2, -0.5) [text=natureBlue, font=\small] {$\mathcal{M}_{\boldsymbol{\theta}}$};

        \coordinate (Center) at (2.2, 1.2);

        \coordinate (P1) at ($(Center)+(-1.8, -0.5)$);
        \coordinate (P2) at ($(Center)+(1.5, -0.5)$);
        \coordinate (P3) at ($(Center)+(1.8, 1.3)$);
        \coordinate (P4) at ($(Center)+(-1.5, 1.3)$);

        \fill[natureGreen!5, opacity=0.6] (P1) -- (P2) -- (P3) -- (P4) -- cycle;
        \draw[thin, natureGreen!50] (P1) -- (P2) -- (P3) -- (P4) -- cycle;
        \node[natureGreen!60!black, font=\scriptsize, anchor=south east] at (P3) {$T_{\hat{\psi}}\mathcal{M}$};

        \coordinate (TrueTip) at ($(Center)+(0.2, 2.5)$);
        \coordinate (NGDTip) at ($(Center)+(1.6, 1.0)$);
        \coordinate (StdTip) at ($(Center)+(1.6, -0.6)$);

        \draw[->, thick, dashed, gray!80, >={Latex[length=2.5mm]}] (Center) -- (StdTip)
        node[pos=1.0, right=2pt, font=\scriptsize, align=left, color=gray] {\textbf{Standard}\\\textbf{Gradient}};

        \draw[->, very thick, natureBlue, >={Latex[length=3mm]}] (Center) -- (NGDTip)
        node[pos=1.0, right=1pt, font=\scriptsize, align=left, color=natureBlue] {\textbf{Projected}\\\textbf{(NGD)}};

        \draw[dashed, natureGreen, thick] (TrueTip) -- (NGDTip)
        node[midway, above right=-2pt, font=\scriptsize] {$\mathcal{P}$};

        \draw[natureGreen] ($(NGDTip)!0.15!(Center)$) -- ++(-0.06, 0.12) -- ($(NGDTip)!0.18!(TrueTip)$);

        \draw[->, very thick, natureRed, >={Latex[length=3mm]}] (Center) -- (TrueTip)
        node[midway, left=2pt, font=\scriptsize, align=right] {\textbf{Sobolev}\\\textbf{Flow}};

        \fill[black] (Center) circle (2.0pt) node[below left=1pt, font=\small] {$\hat{\psi}_{\boldsymbol{\theta}^k}$};

      \end{scope}
    \end{tikzpicture}
  }
  \caption{Schematic of Projected Sobolev NGD method. (a) Commutative diagram relating the continuous function-space flow to the discrete parameter update. (b) Galerkin projection of the Sobolev flow onto the neural tangent space.}
  \label{fig:PSNGD}
\end{figure}

The derivation above holds for a general inner product $(\cdot,\cdot)_X$. In this work, we choose the energy-adaptive Sobolev inner product from \eqref{eq:a_inner}; namely, for $\psi=\hat\psi_{\boldsymbol{\theta}}$, we set $(u,v)_X=a_{\hat\psi_{\boldsymbol{\theta}}}(u,v)$ in the Gram matrix \eqref{eq:gram_matrix_revised}. By contrast, the standard $L^2$ metric yields the classical Fisher information matrix and recovers the stochastic reconfiguration method~\cite{sorella1998green}, but it often induces numerical stiffness that necessitates Runge--Kutta schemes or specialized integrators~\cite{du2021evolutional,bruna2024neural}. The energy-adaptive metric mitigates this stiffness and permits stable explicit integration with large time steps of order $\mathcal{O}(1)$, cf.~\eqref{eq:tau_k}. The subsequent section details efficient strategies for preconditioning the inverse metric $\mathbf{G}(\boldsymbol{\theta}^k)^{-1}$ and the gradient $\nabla_{\boldsymbol{\theta}} L(\boldsymbol{\theta}^k)$.

\section{Implementation via variational Monte Carlo}
\label{sec:vmc}

This section describes the practical implementation of the method: evaluating the loss $L(\boldsymbol{\theta})=E(\hat{\psi}_{\boldsymbol{\theta}})$, estimating its gradient $\nabla_{\boldsymbol{\theta}} L(\boldsymbol{\theta})$, and solving the Gram system involving $\mathbf{G}(\boldsymbol{\theta})$. We first estimate the energy and gradient within a VMC framework using MCMC sampling. We then introduce a hybrid sampling strategy for the normalization constant. Finally, we describe matrix-free Gram products and Nystr\"om preconditioning.

\subsection{Estimation of the energy functional}
\label{subsec:energy-estimation}

Recall that the loss function is defined as the Gross--Pitaevskii energy evaluated at the normalized neural wavefunction $\hat{\psi}_{\boldsymbol{\theta}}$, namely
\begin{equation}
  L(\boldsymbol{\theta}) = E(\hat{\psi}_{\boldsymbol{\theta}})
  =
  \int_{\calD}
  \left(
    \frac{1}{2} |\nabla \hat{\psi}_{\boldsymbol{\theta}}|^2
    +
    V(\mathbf{r}) |\hat{\psi}_{\boldsymbol{\theta}}|^2
    +
    \frac{\beta}{2} |\hat{\psi}_{\boldsymbol{\theta}}|^4
  \right)
  \mathrm{d}\mathbf{r}.
\end{equation}

Let us define the normalization constant (or partition function)
\begin{equation}
  Z(\boldsymbol{\theta})
  =
  \|\psi_{\boldsymbol{\theta}}\|_{L^2(\mathcal D)}^2
  =
  \int_{\calD} |\psi_{\boldsymbol{\theta}}(\mathbf{r})|^2 \, \mathrm{d}\mathbf{r}.
\end{equation}
Within the VMC framework, integrals over $\calD$ are reformulated as expectations with respect to the density $p_{\boldsymbol{\theta}}(\mathbf{r}) = |\hat{\psi}_{\boldsymbol{\theta}}(\mathbf{r})|^2$.

We decompose the total energy into linear and nonlinear components. The linear energy $E_{\mathrm{lin}}$ encompasses the kinetic and external potential terms. By introducing the local energy $E_L(\mathbf{r};\boldsymbol{\theta}) = \psi_{\boldsymbol{\theta}}^{-1} (-\tfrac{1}{2}\nabla^2 + V) \psi_{\boldsymbol{\theta}}$, it takes the expectation form
\begin{equation}
  E_{\mathrm{lin}}
  =
  \int_{\mathcal{D}} \hat{\psi}_{\boldsymbol{\theta}} \left( -\frac{1}{2}\nabla^2 + V \right) \hat{\psi}_{\boldsymbol{\theta}} \, \mathrm{d}\mathbf{r}
  =
  \mathbb{E}_{\mathbf{r}\sim p_{\boldsymbol{\theta}}}
  \bigl[ E_L(\mathbf{r};\boldsymbol{\theta}) \bigr].
\end{equation}
Here $E_L$ is understood pointwise as the ratio $E_L(\mathbf{r};\boldsymbol{\theta}) = \bigl[(-\tfrac{1}{2}\nabla^2 + V)\psi_{\boldsymbol{\theta}}\bigr](\mathbf{r})/\psi_{\boldsymbol{\theta}}(\mathbf{r})$. This ratio is evaluated only away from the nodes; since $p_{\boldsymbol{\theta}}=|\hat\psi_{\boldsymbol{\theta}}|^2$ assigns zero probability to the set $\{\mathbf{r}\in\calD:\psi_{\boldsymbol{\theta}}(\mathbf{r})=0\}$, sampled points avoid the nodes $p_{\boldsymbol{\theta}}$-almost surely.
Notably, $E_{\mathrm{lin}}$ is invariant under the rescaling of $\psi_{\boldsymbol{\theta}}$, so $Z(\boldsymbol{\theta})$ cancels out. In contrast, the nonlinear interaction energy $E_{\mathrm{int}}$ explicitly depends on the amplitude
\begin{equation}
  E_{\mathrm{int}}
  =
  \frac{\beta}{2}
  \int_{\mathcal{D}} |\hat{\psi}_{\boldsymbol{\theta}}|^4 \,\mathrm{d}\mathbf{r}
  =
  \mathbb{E}_{\mathbf{r}\sim p_{\boldsymbol{\theta}}}
  \left[
    \frac{\beta}{2 Z(\boldsymbol{\theta})}
    |\psi_{\boldsymbol{\theta}}(\mathbf{r})|^2
  \right].
\end{equation}
Combining these contributions yields the total energy expectation
\begin{equation}
  L(\boldsymbol{\theta})
  =
  \mathbb{E}_{\mathbf{r}\sim p_{\boldsymbol{\theta}}}
  \left[
    E_L(\mathbf{r};\boldsymbol{\theta})
    +
    \frac{\beta}{2 Z(\boldsymbol{\theta})}
    |\psi_{\boldsymbol{\theta}}(\mathbf{r})|^2
  \right].
  \label{eq:total_energy_expectation}
\end{equation}
This formulation highlights that a stable evaluation of the scalar $Z(\boldsymbol{\theta})$ is essential. In \Cref{sec:ADIS}, we describe a hybrid sampling strategy designed to compute $Z(\boldsymbol{\theta})$ effectively.

\subsection{Estimation of the gradients}\label{subsec:grad}

To estimate the gradient $\nabla_{\boldsymbol{\theta}} L(\boldsymbol{\theta})$, we employ the log-derivative identity $\nabla_{\boldsymbol{\theta}} |\psi_{\boldsymbol{\theta}}|^2 = 2 |\psi_{\boldsymbol{\theta}}|^2 \mathbf{O}_{\boldsymbol{\theta}}$, where $\mathbf{O}_{\boldsymbol{\theta}}(\mathbf{r}) = \nabla_{\boldsymbol{\theta}} \ln |\psi_{\boldsymbol{\theta}}(\mathbf{r})|$ denotes the score function. Differentiating the linear Rayleigh quotient $E_{\mathrm{lin}}$ and the nonlinear interaction energy $E_{\mathrm{int}}$ with respect to $\boldsymbol{\theta}$, and applying the quotient rule to handle the normalization factor $Z(\boldsymbol{\theta})$, allows us to express the total gradient as an expectation over the probability density $p_{\boldsymbol{\theta}}$.

To present the estimator in a concise form, we introduce the \emph{local chemical potential} $\mu_L(\mathbf{r})$, which effectively aggregates the local energy and the potential induced by the nonlinearity
\begin{equation}
  \mu_L(\mathbf{r}) \coloneqq E_L(\mathbf{r}) + \frac{\beta}{Z(\boldsymbol{\theta})} |\psi_{\boldsymbol{\theta}}(\mathbf{r})|^2.
\end{equation}
Correspondingly, we denote the mean chemical potential as $\bar{\mu} = \mathbb{E}_{\mathbf{r} \sim p_{\boldsymbol{\theta}}}[\mu_L(\mathbf{r})]$, noting that $\bar{\mu} = E_{\mathrm{lin}} + 2 E_{\mathrm{int}}$ due to the scaling of the interaction term. Combining these definitions, the gradient of the loss function is given by
\begin{equation}
  \nabla_{\boldsymbol{\theta}} L(\boldsymbol{\theta})
  =
  2 \,
  \mathbb{E}_{\mathbf{r} \sim p_{\boldsymbol{\theta}}}
  \left[
    (\mu_L(\mathbf{r}) - \bar{\mu})
    \, \mathbf{O}_{\boldsymbol{\theta}}(\mathbf{r})
  \right].
  \label{eq:grad_total_final}
\end{equation}
This estimator, structured as a covariance between the local chemical potential and the score function, inherently acts as a control variate formulation, reducing the variance of the stochastic gradients.

\subsection{Hybrid sampling strategy}\label{sec:ADIS}

The explicit dependence of the energy expectation \eqref{eq:total_energy_expectation} on the normalization constant $Z(\boldsymbol{\theta})$ necessitates a dual sampling approach. We employ a hybrid strategy comprising two distinct routines:

\textbf{Gradient Estimation:} The first routine estimates energy gradients via MCMC sampling (e.g., Metropolis--Hastings (M-H) \cite{metropolis1953equation} or Hamiltonian Monte Carlo (HMC) \cite{neal2011mcmc}) from $p_{\boldsymbol{\theta}}$. By concentrating samples in high-probability regions where $|\hat{\psi}_{\boldsymbol{\theta}}|^2$ is large, this method yields low-variance estimators efficiently.

\textbf{Normalization via ADIS:} The second routine computes $Z(\boldsymbol{\theta})$ using Monte Carlo integration to overcome the curse of dimensionality inherent in grid-based quadrature. We utilize an adaptive defensive importance sampling (ADIS) strategy, expressing the constant as $Z(\boldsymbol{\theta}) = \mathbb{E}_{\mathbf{r} \sim q} [ |\psi_{\boldsymbol{\theta}}|^2 / q ]$. The proposal $q$ is a mixture
\begin{equation}
  q(\mathbf{r}) = \alpha \, q_{\mathcal{N}}(\mathbf{r} \mid \boldsymbol{\mu}, \boldsymbol{\Sigma}) + (1-\alpha) \, q_{\mathcal{U}}(\mathbf{r} \mid \Omega), \quad \alpha \in [0,1].
\end{equation}
Here, $q_{\mathcal{N}}$ is a multivariate Gaussian updated using statistics from the MCMC stream to track the wavefunction's mass, while the defensive uniform component $q_{\mathcal{U}}$ guarantees global support. This combination prevents variance explosion and ensures the method scales linearly with sample size, independent of spatial dimension.
All reported main experiments use the fixed default $\alpha=0.8$, with no per-problem tuning. The ablations in \Cref{sm:adis-alpha} show that this value gives the most accurate $Z$ estimate in the 3D benchmark with a grid reference and remains stable in the $d=8$ benchmark.

\subsection{Computation of the Gram matrix}

We adopt the $a_{\hat{\psi}}$-inner product evaluated at the normalized wavefunction $\hat{\psi}_{\boldsymbol{\theta}}$ as the underlying metric for constructing the Gram matrix, cf.~\eqref{eq:a_inner}. Specifically, the bilinear form is defined as
\begin{equation}
  a_{\hat{\psi}}(u,v)
  = \int_{\calD} \left(
    \frac{1}{2}\nabla u \cdot \nabla v
    + V(\mathbf{r})\, u v
    + \frac{\beta}{Z(\boldsymbol{\theta})} |\psi_{\boldsymbol{\theta}}(\mathbf{r})|^2 u v
  \right)\,\mathrm{d}\mathbf{r},
  \label{eq:metric_def_ref}
\end{equation}
The Gram matrix is then given by
\begin{equation}
  \mathbf{G}_{ij}(\boldsymbol{\theta})
  = a_{\hat{\psi}}\!\left(
    \frac{\partial \hat{\psi}_{\boldsymbol{\theta}}}{\partial \theta_i},
    \frac{\partial \hat{\psi}_{\boldsymbol{\theta}}}{\partial \theta_j}
  \right).
\end{equation}

We first derive an explicit expression for the tangent vectors associated with the normalized wavefunction $\hat{\psi}_{\boldsymbol{\theta}} = \psi_{\boldsymbol{\theta}} / \sqrt{Z(\boldsymbol{\theta})}$. Differentiating with respect to a parameter $\theta_k$ yields
\begin{align}
  \frac{\partial \hat{\psi}_{\boldsymbol{\theta}}}{\partial \theta_k}
  &= \frac{\partial}{\partial \theta_k}\!\left(\frac{\psi_{\boldsymbol{\theta}}}{\sqrt{Z}}\right)
  = \frac{(\partial_{\theta_k}\psi_{\boldsymbol{\theta}})\sqrt{Z}
  - \psi_{\boldsymbol{\theta}} \left(\frac{1}{2\sqrt{Z}}\partial_{\theta_k} Z\right)}{Z}
  \nonumber \\
  &= \hat{\psi}_{\boldsymbol{\theta}}
  \left(
    \frac{\partial_{\theta_k}\psi_{\boldsymbol{\theta}}}{\psi_{\boldsymbol{\theta}}}
    - \frac{1}{2}\frac{\partial_{\theta_k} Z}{Z}
  \right).
\end{align}
Introducing the log-derivative $
\mathbf{O}_k(\mathbf{r}) := \partial_{\theta_k}\ln|\psi_{\boldsymbol{\theta}}(\mathbf{r})|,
$
and recalling that $\partial_{\theta_k} Z = 2Z\,\mathbb{E}_{\mathbf{r}\sim p_{\boldsymbol{\theta}}}[\mathbf{O}_k]$, we arrive at the compact form
\begin{equation}
  \frac{\partial \hat{\psi}_{\boldsymbol{\theta}}}{\partial \theta_k}
  = \hat{\psi}_{\boldsymbol{\theta}}(\mathbf{r})
  \left(
    \mathbf{O}_k(\mathbf{r})
    - \mathbb{E}_{\mathbf{r}\sim p_{\boldsymbol{\theta}}}[\mathbf{O}_k]
  \right)
  =: \hat{\psi}_{\boldsymbol{\theta}}(\mathbf{r})\,\delta \mathbf{O}_k(\mathbf{r}),
  \label{eq:tangent_vec_centered}
\end{equation}
where $\delta \mathbf{O}_k$ denotes the \emph{centered score function}, satisfying $\mathbb{E}_{p_{\boldsymbol{\theta}}}[\delta \mathbf{O}_k]=0$. This centering reflects the orthogonality of tangent vectors to the normalization constraint.

\subsubsection{Gradient expansion of the kinetic term}

We now evaluate the kinetic contribution
\[
  \frac{1}{2}\nabla(\partial_{\theta_i}\hat{\psi}_{\boldsymbol{\theta}})\cdot
  \nabla(\partial_{\theta_j}\hat{\psi}_{\boldsymbol{\theta}}).
\]
Using \eqref{eq:tangent_vec_centered} and the product rule,
\begin{equation}
  \nabla\!\left(\hat{\psi}_{\boldsymbol{\theta}} \,\delta \mathbf{O}_k\right)
  = (\nabla\hat{\psi}_{\boldsymbol{\theta}})\,\delta \mathbf{O}_k
  + \hat{\psi}_{\boldsymbol{\theta}}\,\nabla\delta \mathbf{O}_k.
\end{equation}
Define the \emph{quantum force} $
\mathbf{F}(\mathbf{r})
:= \nabla \ln|\hat{\psi}_{\boldsymbol{\theta}}(\mathbf{r})|
= \frac{\nabla\hat{\psi}_{\boldsymbol{\theta}}}{\hat{\psi}_{\boldsymbol{\theta}}}.
$
Since $\delta \mathbf{O}_k = \mathbf{O}_k - \mathbb{E}[\mathbf{O}_k]$, its gradient satisfies $\nabla\delta \mathbf{O}_k = \nabla \mathbf{O}_k$. Therefore,
\begin{equation}
  \nabla\!\left(\frac{\partial \hat{\psi}_{\boldsymbol{\theta}}}{\partial \theta_k}\right)
  = \hat{\psi}_{\boldsymbol{\theta}}(\mathbf{r})
  \left(
    \mathbf{F}(\mathbf{r})\,\delta \mathbf{O}_k(\mathbf{r})
    + \nabla \mathbf{O}_k(\mathbf{r})
  \right).
\end{equation}
Substituting the expressions for indices $i$ and $j$ into the kinetic term yields
\begin{align}
  \mathbf{K}_{ij}
  &= \int_{\calD}
  \frac{1}{2}
  \left[\hat{\psi}_{\boldsymbol{\theta}}(\mathbf{F}\delta \mathbf{O}_i + \nabla \mathbf{O}_i)\right]
  \cdot
  \left[\hat{\psi}_{\boldsymbol{\theta}}(\mathbf{F}\delta \mathbf{O}_j + \nabla \mathbf{O}_j)\right]
  \,\mathrm{d}\mathbf{r}
  \nonumber \\
  &= \int_{\calD}
  \frac{1}{2}|\hat{\psi}_{\boldsymbol{\theta}}|^2
  \Big[
    |\mathbf{F}|^2\,\delta \mathbf{O}_i\delta \mathbf{O}_j
    + \mathbf{F}\cdot(\delta \mathbf{O}_i\nabla \mathbf{O}_j + \delta \mathbf{O}_j\nabla \mathbf{O}_i)
    + \nabla \mathbf{O}_i\cdot\nabla \mathbf{O}_j
  \Big]
  \,\mathrm{d}\mathbf{r}
  \nonumber \\
  &= \mathbb{E}_{\mathbf{r}\sim p_{\boldsymbol{\theta}}}\!\left[
    \frac{1}{2}|\mathbf{F}|^2\,\delta \mathbf{O}_i\delta \mathbf{O}_j
    + \frac{1}{2}\mathbf{F}\cdot(\delta \mathbf{O}_i\nabla \mathbf{O}_j + \delta \mathbf{O}_j\nabla \mathbf{O}_i)
    + \frac{1}{2}\nabla \mathbf{O}_i\cdot\nabla \mathbf{O}_j
  \right],
  \label{eq:k}
\end{align}
where we have used $p_{\boldsymbol{\theta}} = |\hat{\psi}_{\boldsymbol{\theta}}|^2$ to convert the integral into an expectation.

\subsubsection{Potential and interaction terms}
Substituting the tangent vectors $\hat{\psi}_{\boldsymbol{\theta}}\,\delta \mathbf{O}_k$ directly, we obtain
\begin{align}
  \mathbf{V}_{ij} + \mathbf{I}_{ij}
  &= \int_{\calD}
  \left(
    V(\mathbf{r})
    + \frac{\beta}{Z(\boldsymbol{\theta})}|\psi_{\boldsymbol{\theta}}(\mathbf{r})|^2
  \right)
  (\hat{\psi}_{\boldsymbol{\theta}}\,\delta \mathbf{O}_i)
  (\hat{\psi}_{\boldsymbol{\theta}}\,\delta \mathbf{O}_j)
  \,\mathrm{d}\mathbf{r}
  \nonumber \\
  &= \mathbb{E}_{\mathbf{r}\sim p_{\boldsymbol{\theta}}}\!\left[
    \left(
      V(\mathbf{r})
      + \frac{\beta}{Z(\boldsymbol{\theta})}|\psi_{\boldsymbol{\theta}}(\mathbf{r})|^2
    \right)
    \delta \mathbf{O}_i(\mathbf{r})\,\delta \mathbf{O}_j(\mathbf{r})
  \right].
  \label{eq:vandi}
\end{align}

Collecting \eqref{eq:k} and \eqref{eq:vandi}, it is convenient to group all coefficients multiplying $\delta \mathbf{O}_i\delta \mathbf{O}_j$ into a single weight function, termed the \emph{generalized energy density}
\begin{equation}
  \mathcal{W}(\mathbf{r})
  =V(\mathbf{r})
  + \frac{\beta}{Z(\boldsymbol{\theta})}|\psi_{\boldsymbol{\theta}}(\mathbf{r})|^2
  +\frac{1}{2}|\mathbf{F}(\mathbf{r})|^2.
\end{equation}
The Sobolev Gram matrix then admits the compact expectation form
\begin{equation}
  \label{eq:sobolev_gram_final}
  \begin{aligned}
    \mathbf{G}_{ij}(\boldsymbol{\theta})
    = \mathbb{E}_{\mathbf{r}\sim p_{\boldsymbol{\theta}}}\Big[
      & \underbrace{\mathcal{W}(\mathbf{r})\,\delta \mathbf{O}_i(\mathbf{r})\,\delta \mathbf{O}_j(\mathbf{r})}_{\text{Term I: energy-weighted Fisher information}} \\
      & + \underbrace{\frac{1}{2}\mathbf{F}(\mathbf{r})\cdot
        \big(\nabla \mathbf{O}_i(\mathbf{r})\,\delta \mathbf{O}_j(\mathbf{r})
      + \nabla \mathbf{O}_j(\mathbf{r})\,\delta \mathbf{O}_i(\mathbf{r})\big)}_{\text{Term II: drift--gradient coupling}} \\
      & + \underbrace{\frac{1}{2}\nabla \mathbf{O}_i(\mathbf{r})\cdot\nabla \mathbf{O}_j(\mathbf{r})}_{\text{Term III: pure Sobolev regularization}}
    \Big].
  \end{aligned}
\end{equation}

\subsection{Matrix-free PCG solver with Nyström preconditioning}
\label{subsec:nystrom}

The Gram matrix $\mathbf{G}(\boldsymbol{\theta})$ is symmetric positive semidefinite. To make the linear system nonsingular, we add Tikhonov regularization $\lambda \mathbf{I}$ with $\lambda>0$. Projecting the Sobolev gradient flow onto parameter space then requires solving $(\mathbf{G}(\boldsymbol{\theta}) + \lambda \mathbf{I})\,\dot{\boldsymbol{\theta}} = - \nabla_{\boldsymbol{\theta}} L(\boldsymbol{\theta})$. Explicitly forming the Gram matrix costs $\mathcal{O}(P^2)$ memory, and direct inversion costs $\mathcal{O}(P^3)$ time. We therefore use a matrix-free PCG solver that requires only products $\mathbf{v} \mapsto \mathbf{G}\mathbf{v}$.

\subsubsection{Matrix-free implementation}
\label{subsubsec:matrix-free-implementation}

For Monte Carlo samples $\{\mathbf{r}_n\}_{n=1}^N$, let
\[
\mathbf{O}(\mathbf{r})=\nabla_{\boldsymbol{\theta}}\log|\psi_{\boldsymbol{\theta}}(\mathbf{r})|\in\mathbb{R}^P,
\qquad
\delta\mathbf{O}(\mathbf{r})=\mathbf{O}(\mathbf{r})-\mathbb{E}_{p_{\boldsymbol{\theta}}}[\mathbf{O}]
\]
be the score vector and its centered version. Define the centered score matrix and mixed input--parameter derivative tensor by
\[
\mathbf{J}\in\mathbb{R}^{N\times P},
\qquad
\mathbf{J}_{n,:}=\delta\mathbf{O}(\mathbf{r}_n)^\top,
\]
\[
\mathbf{K}\in\mathbb{R}^{N\times P\times D},
\qquad
\mathbf{K}_{n,:,d}
=\nabla_{r_d}\delta\mathbf{O}(\mathbf{r}_n)^\top
=\nabla_{r_d}\mathbf{O}(\mathbf{r}_n)^\top .
\]
The last equality holds because the centering constant has no spatial derivative. Let $\mathbf{F}(\mathbf{r}_n)\in\mathbb{R}^D$ denote the drift field, and let $\mathcal{W}(\mathbf{r}_n)$ denote the generalized energy density. The Sobolev Gram matrix in \eqref{eq:sobolev_gram_final} acts on $\mathbf{v}\in\mathbb{R}^P$ as
\begin{equation}
\label{eq:mvp-explicit-contraction}
  \begin{aligned}
    \mathbf{G}\mathbf{v}
    =
    \frac{1}{N}
    \Big[
      & \mathbf{J}^\top \bigl( \mathbf{W} \odot (\mathbf{J} \mathbf{v}) \bigr)
      + \frac12 \sum_{d=1}^D
      \Big(
        \mathbf{K}_{:,:,d}^\top
        \bigl(
          \mathbf{F}_{:,d} \odot (\mathbf{J} \mathbf{v})
        \bigr)
        +
        \mathbf{J}^\top
        \bigl(
          \mathbf{F}_{:,d} \odot (\mathbf{K}_{:,:,d} \mathbf{v})
        \bigr)
      \Big)
      \\[0.3em]
      & + \frac12 \sum_{d=1}^D
      \mathbf{K}_{:,:,d}^\top (\mathbf{K}_{:,:,d} \mathbf{v})
    \Big],
  \end{aligned}
\end{equation}
where $\mathbf{W}\in\mathbb{R}^N$ stacks the generalized energy densities and $\odot$ denotes elementwise multiplication. We now describe two implementations of this product, which differ in whether the matrix $\mathbf{J}$ and the mixed-derivative tensor $\mathbf{K}$ are materialized.

\textit{Explicit Sobolev implementation.} In the low-dimensional experiments, $\mathbf{J}$ and $\mathbf{K}$ are materialized and the Gram products are evaluated by batched contractions \eqref{eq:mvp-explicit-contraction}. The centered score $\delta\mathbf{O}(\mathbf{r}_n)$ is obtained by a reverse-mode pass over the parameters, with cost $\mathcal{O}(P)$ per sample. The tensor slice $\mathbf{K}_{n,:,d}=\nabla_{r_d}\mathbf{O}(\mathbf{r}_n)^\top$ is computed by forward-over-reverse automatic differentiation. Consequently, forming $\mathbf{K}$ costs $\mathcal{O}(NDP)$ operations, up to constants from automatic differentiation, and requires $\mathcal{O}(NPD)$ storage. Once $\mathbf{J}$ and $\mathbf{K}$ are cached, each product $\mathbf{v}\mapsto\mathbf{G}\mathbf{v}$ costs $\mathcal{O}(NDP)$ floating-point operations. This path avoids explicit $\mathcal{O}(P^2)$ Gram storage and $\mathcal{O}(P^3)$ Gram inversion, but storing $\mathbf{K}$ becomes prohibitive as $P$ or $D$ grows.

\textit{Implicit Sobolev operator chain.} The implicit implementation evaluates the same product without materializing either $\mathbf{J}$ or $\mathbf{K}$.
 It uses the quadratic-form identity
\[
\mathbf{G}\mathbf{v}
=
\nabla_{\mathbf{v}}\left[\frac12\mathbf{v}^\top\mathbf{G}\mathbf{v}\right].
\]
Given $\mathbf{v}\in\mathbb{R}^P$, a forward-mode Jacobian-vector product first computes the uncentered directional score
\[
O_{\mathbf{v}}(\mathbf{r}_n)
=
\sum_j v_j\,\partial_{\theta_j}\log|\psi_{\boldsymbol{\theta}}(\mathbf{r}_n)| .
\]
Sample-mean centering gives the centered tangent perturbation
\[
  \delta\psi_{\mathbf{v}}(\mathbf{r}_n)
  =
  O_{\mathbf{v}}(\mathbf{r}_n)
  -
  \frac1N\sum_{m=1}^N O_{\mathbf{v}}(\mathbf{r}_m)
  =
  \sum_j v_j\,\delta\mathbf{O}_j(\mathbf{r}_n).
\]
Spatial differentiation then yields the input-gradient intermediate $\nabla_{\mathbf{r}}\delta\psi_{\mathbf{v}}\in\mathbb{R}^{N\times D}$. Writing $\delta\psi_{\mathbf{v},n}=\delta\psi_{\mathbf{v}}(\mathbf{r}_n)$, these quantities define the scalar tangent energy
\begin{equation}
  \label{eq:tangent-energy}
  \mathcal{E}(\mathbf{v})
  =
  \frac1N\sum_{n=1}^N\Bigl[
    \tfrac12\mathcal{W}(\mathbf{r}_n)\,\delta\psi_{\mathbf{v},n}^2
    +\tfrac12\mathbf{F}(\mathbf{r}_n)\cdot\nabla_{\mathbf{r}}\delta\psi_{\mathbf{v},n}\,\delta\psi_{\mathbf{v},n}
    +\tfrac14\|\nabla_{\mathbf{r}}\delta\psi_{\mathbf{v},n}\|^2
  \Bigr]
  =
  \tfrac12\mathbf{v}^\top\mathbf{G}\mathbf{v} .
\end{equation}
A final reverse-mode vector-Jacobian product returns $\mathbf{G}\mathbf{v}=\nabla_{\mathbf{v}}\mathcal{E}(\mathbf{v})$.   The product rule reproduces the two drift--gradient cross terms in \eqref{eq:mvp-explicit-contraction}. This is Pearlmutter's Hessian-vector trick \cite{pearlmutter1994fast} specialized to the energy-adaptive metric. Each matrix-free product stores only samplewise directional intermediates such as $O_{\mathbf v}\in\mathbb R^N$ and $\nabla_{\mathbf r}O_{\mathbf v}\in\mathbb R^{N\times D}$, parameter-sized gradient vectors, the Nystr\"om basis $\mathbf M\in\mathbb{R}^{P\times r}$, and automatic-differentiation activations. Its peak-memory scaling is therefore $\mathcal{O}(N)+\mathcal{O}(ND)+\mathcal{O}(Pr)+\text{activations}$, with no $\mathcal O(NP)$ or $\mathcal O(NPD)$ cached tensor. It also avoids the dense $\mathcal{O}(P^2)$ Gram matrix.

The standard Fisher or stochastic-reconfiguration metric uses only the first-order centered score Jacobian $\mathbf{J}$, so its matrix-free product costs $\mathcal{O}(NP)$. In our implementation, one implicit Sobolev product empirically corresponds to roughly four to six network forward/backward evaluations. The constant $c_{\mathrm{AD}}$ in \Cref{tab:mvp_complexity} accounts for the JVP, input-gradient, and VJP chain. We use the implicit path in the high-dimensional experiments.

\begin{table}[!htbp]
  \centering
  \caption{Cost of one matrix-free product $\mathbf{v}\mapsto\mathbf{G}\mathbf{v}$ and peak memory of the associated Nystr\"om--PCG solve. Here $N$ is the number of samples, $P$ the number of parameters, $D$ the spatial dimension, $r$ the Nystr\"om rank, and $c_{\mathrm{AD}}$ an implementation-dependent automatic-differentiation constant.}
  \label{tab:mvp_complexity}
  \resizebox{\textwidth}{!}{
  \begin{tabular}{lcc}
    \toprule
     & Explicit Sobolev & Implicit Sobolev \\
    \midrule
    Jacobian $\mathbf{J}$        & $\mathcal{O}(NP)$, materialized  & not materialized \\
    Mixed derivative $\mathbf{K}$ & $\mathcal{O}(NPD)$, materialized & not materialized \\
    Product mechanism & tensor contractions & JVP/input gradient/VJP \\
    Per-product work     & $\mathcal{O}(NDP)$ & $\mathcal{O}(NP\,c_{\mathrm{AD}})$ \\
    Peak memory                  & $\mathcal{O}(NPD)+\mathcal{O}(Pr)$ + activations & $\mathcal{O}(N)+\mathcal{O}(ND)+\mathcal{O}(Pr)$ + activations \\
    \bottomrule
  \end{tabular}}
\end{table}

\subsubsection{Randomized Nyström preconditioning}
\label{subsubsec:nystrom-preconditioning}

To accelerate PCG convergence, we use a low-rank Nyström approximation of the Gram matrix \cite{frangella2023randomized}: $\mathbf{G} \approx \mathbf{M} \mathbf{M}^\top$ with $\mathbf{M} \in \mathbb{R}^{P \times r}$ and $r \ll P$. The construction follows the GPU-efficient approach of \cite{guzman2025improving}:
\begin{enumerate}
  \item Generate a Gaussian random test matrix $\boldsymbol{\Omega} \in \mathbb{R}^{P \times r}$;
  \item Compute the projected matrix $\mathbf{Y} = \mathbf{G} \boldsymbol{\Omega}$ using the matrix-free $\mathbf{G}\mathbf{v}$ operator;
  \item Form the Nyström basis
    \[
      \mathbf{M} = \mathbf{Y} \bigl(\boldsymbol{\Omega}^\top \mathbf{Y} + \varepsilon \mathbf{I}\bigr)^{-1/2},
    \]
    where a small shift $\varepsilon>0$ is added for numerical stability.
\end{enumerate}

The Woodbury identity gives a closed form for the approximate inverse used as the preconditioner:
\begin{equation}
  (\mathbf{M} \mathbf{M}^\top + \lambda \mathbf{I})^{-1}
  = \frac{1}{\lambda}
  \left[
    \mathbf{I} - \mathbf{M} \left( \mathbf{M}^\top \mathbf{M} + \lambda \mathbf{I}_r \right)^{-1} \mathbf{M}^\top
  \right].
\end{equation}
This reduces the inversion cost to $\mathcal{O}(Pr^2)$, which is tractable for moderate ranks. We also use a lazy update strategy that refreshes the Nyström basis only every ten steps, reducing overhead. The reported main experiments use the fixed default $r=128$ and update frequency $K=10$. The rank ablation in \Cref{sm:nystrom-rank} shows that $r=128$ is a budget-aware default, while $r=512$ gives the best median residual and Hutchinson spectral error for the tested $P\sim 10^4$ network when accuracy is prioritized.

\subsection{Implementation details}
\label{sec:other_tech}

\Cref{alg:sobolev_ngd} summarizes the Sobolev NGD algorithm. The implementation combines several techniques to improve robustness and scalability.
\begin{algorithm}[!htbp]
  \caption{Projected Sobolev NGD for the GPE}
  \label{alg:sobolev_ngd}
  \begin{algorithmic}[1]
    \Require Initial parameters $\boldsymbol{\theta}_0$, learning rate $\tau_k$, damping parameter $\lambda_k$, Nyström rank $r=128$, and update frequency $K=10$.
    \Ensure Ground-state wavefunction approximation $\hat{\psi}_{\boldsymbol{\theta}}$.

    \State \textbf{Initialization}: Initialize MCMC walkers $\mathcal{X}_{\text{mcmc}}$.
    \For{$k = 0, 1, 2, \dots$ until convergence}
    \State \textbf{Step 1: Hybrid Sampling and Statistics}
    \State \quad Update walkers $\mathcal{X}_{\text{mcmc}} = \{\mathbf{r}_n\}_{n=1}^N$ via MCMC sampling from $p_{\boldsymbol{\theta}}(\mathbf{r}) \propto |\psi_{\boldsymbol{\theta}}(\mathbf{r})|^2$.
    \State \quad Compute statistics $(\boldsymbol{\mu}, \boldsymbol{\Sigma})$ from $\mathcal{X}_{\text{mcmc}}$ and update the ADIS proposal $q(\mathbf{r})$.
    \State \quad Estimate the normalization constant $Z(\boldsymbol{\theta})$ using ADIS samples $\mathcal{X}_{\text{int}}$.

    \State \textbf{Step 2: Gradient Estimation}
    \State \quad Evaluate local chemical potential $\mu_L(\mathbf{r}_n) = E_L(\mathbf{r}_n) + \frac{\beta}{Z(\boldsymbol{\theta})} |\psi_{\boldsymbol{\theta}}(\mathbf{r}_n)|^2$.
    \State \quad Form the Monte Carlo loss-gradient estimate $\widehat{\mathbf{g}}_k = \frac{2}{N} \sum_{n=1}^N (\mu_L(\mathbf{r}_n) - \bar{\mu}) \nabla_{\boldsymbol{\theta}} \ln |\psi_{\boldsymbol{\theta}}(\mathbf{r}_n)|$.

    \State \textbf{Step 3: Randomized Nyström Preconditioner (Lazy Update)}
    \If{$k \equiv 0 \pmod K$}
    \State Generate Gaussian random matrix $\boldsymbol{\Omega} \in \mathbb{R}^{P \times r}$.
    \State Compute $\mathbf{Y} = \mathbf{G} \boldsymbol{\Omega}$ via the matrix-free operator.
    \State Form the Nyström basis $\mathbf{M} = \mathbf{Y} (\boldsymbol{\Omega}^\top \mathbf{Y} + \varepsilon \mathbf{I})^{-1/2}$.
    \EndIf

    \State \textbf{Step 4: Natural Gradient Direction Search}
    \State \quad Solve $(\mathbf{G}_k + \lambda_k \mathbf{I}) \dot{\boldsymbol{\theta}}_k = -\widehat{\mathbf{g}}_k$ using matrix-free PCG.
    \State \quad Apply preconditioning via the Woodbury identity:
    \State \quad $\mathbf{P}^{-1} \mathbf{v} = \frac{1}{\lambda_k} \left[ \mathbf{v} - \mathbf{M} \left( \mathbf{M}^\top \mathbf{M} + \lambda_k \mathbf{I}_r \right)^{-1} \mathbf{M}^\top \mathbf{v} \right]$.

    \State \textbf{Step 5: Parameter and Damping Update}
    \State \quad Update parameters: $\boldsymbol{\theta}_{k+1} = \boldsymbol{\theta}_k + \tau_k \dot{\boldsymbol{\theta}}_k$.
    \State \quad Update damping by the residual-aware rule \eqref{eq:lambda_k}.
    \EndFor
    \State \Return $\hat{\psi}_{\boldsymbol{\theta}} = \psi_{\boldsymbol{\theta}} / \sqrt{Z(\boldsymbol{\theta})}$
  \end{algorithmic}
\end{algorithm}

We first clarify the role of Tikhonov regularization. The continuous projected Sobolev gradient flow \eqref{eq:projected_gradient_flow_revised} is energy dissipative in function space. However, this monotonicity should not be interpreted as a pathwise monotonicity guarantee for the fully discrete stochastic NGD iterates. In practice, both the gradient and the Gram-matrix action are estimated from Monte Carlo samples, and the regularized Gram system is solved only approximately by PCG. The damping parameter defines $\mathbf{H}_k=\mathbf{G}_k+\lambda_k\mathbf{I}$ and therefore satisfies $\|\mathbf{H}_k^{-1}\|\le 1/\lambda_k$. Hence, the lower bound imposed on $\lambda_k$ limits the amplification of both Monte Carlo noise and PCG residuals. The specific residual-aware schedule and its empirical trajectory are reported in \eqref{eq:lambda_k} and \Cref{subfig:lambda_lattice}, respectively. An idealized conditional-expectation stability estimate under a conditionally unbiased gradient model is provided in \Cref{sm:damping-stability}.

We next describe the Fourier-feature embedding used for oscillatory potentials. Because standard MLPs suffer from spectral bias, we use Gaussian random Fourier features (RFFs) \cite{tancik2020fourier} for potentials with known oscillatory lattice structure. The embedding is
\[
  \gamma(\mathbf{r})
  =
  [\cos(2\pi \mathbf{B}\mathbf{r}),\,
   \sin(2\pi \mathbf{B}\mathbf{r})]\in\mathbb{R}^{2m},
\]
where the rows of $\mathbf{B}\in\mathbb{R}^{m\times d}$ are initialized independently from $\mathcal{N}(0,\sigma^2I_d)$. In the optical-lattice experiments, $\mathbf{B}$ is trainable after this Gaussian initialization and is optimized jointly with the MLP. In the smooth high-dimensional anisotropic quadratic experiments, no RFF layer is used; the raw coordinate $\mathbf{r}$ is passed directly to the MLP. This embedding helps the network resolve high-frequency spatial structures when they are present in the potential.

\begin{table}[!htbp]
  \centering
  \caption{Random Fourier feature configurations used in the experiments.}
  \label{tab:rff_config}
  \resizebox{\textwidth}{!}{
  \begin{tabular}{lcccl}
    \toprule
    Experiment & $d$ & $m$ & $\sigma$ & Rationale \\
    \midrule
    2D optical lattice & 2 & 16 & 1.0 & oscillatory term $\sin^2(\pi x)$ \\
    3D optical lattice & 3 & 16 & 0.25 & oscillatory term $\sin^2(\pi x/4)$ \\
    High-dimensional anisotropic quadratic & 4--8 & -- & -- & no RFFs; smooth low-frequency state \\
    \bottomrule
\end{tabular}}
\end{table}

For the optical-lattice benchmarks, the bandwidth is chosen according to the spatial scale
of the prescribed lattice oscillations. We use $\sigma=1.0$ for the 2D $\sin^2(\pi x)$ lattice and the smaller value $\sigma=0.25$ for the slower-varying 3D $\sin^2(\pi x/4)$ lattice. These choices are summarized in \Cref{tab:rff_config}.

\section{Numerical experiments}
\label{sec:numerical_experiments}

This section evaluates the proposed method on problems ranging from low-dimensional optical lattices with strong nonlinearity to high-dimensional systems with fully coupled anisotropic harmonic confinement. We also compare the method with norm-DNN and test whether the neural-network solution can serve as an effective initial guess for traditional high-precision solvers. These experiments assess accuracy, computational efficiency, and scaling with spatial dimension.

\subsection{Experimental setup and evaluation metrics}
\label{subsec:experimental-setup}

We measure accuracy using the variance of the local chemical potential (the \emph{local variance}), which for normalized wavefunctions equals the squared residual norm of the nonlinear eigenvalue problem
\begin{equation}
  \begin{aligned}
    \text{Res}(\hat{\psi})
    &= \|H[\hat{\psi}]\hat{\psi} - \mu_1 \hat{\psi}\|_2^2
    = \int |\hat{\psi}(\mathbf{r})|^2
    \left| \frac{H[\hat{\psi}]\hat{\psi}(\mathbf{r})}{\hat{\psi}(\mathbf{r})} - \mu_1 \right|^2
    \,\mathrm{d}\mathbf{r} \\
    &= \mathbb{E}_{\mathbf{r} \sim |\hat{\psi}|^2}
    \left[ \bigl(\mu_{\mathrm{L}}(\mathbf{r}) - \mu_1\bigr)^2 \right]
    = \text{Var}_{\mathbf{r} \sim |\hat{\psi}|^2}
    \bigl[\mu_{\mathrm{L}}(\mathbf{r})\bigr].
  \end{aligned}
\end{equation}

For the low-dimensional problems (2D and 3D), we further validate the results against high-accuracy reference solutions computed with the \texttt{DFTK.jl} package \cite{DFTKjcon}, which uses direct energy minimization with plane-wave basis sets.
We report the following relative error metrics:
\begin{align}
  \varepsilon_E &= \frac{|E_{\text{NN}} - E_{\text{ref}}|}{|E_{\text{ref}}|},
  &\text{(relative energy error)}, \\
  \varepsilon_\mu &= \frac{|\mu_1^{\text{NN}} - \mu_{\text{ref}}|}{|\mu_{\text{ref}}|},
  &\text{(relative chemical potential error)}, \\
  \varepsilon_\rho &= \frac{\|\rho_{\text{NN}} - \rho_{\text{ref}}\|_2}{\|\rho_{\text{ref}}\|_2},
  &\text{(relative density $L^2$ error)},
\end{align}
where $\rho(\mathbf{r}) = |\hat{\psi}(\mathbf{r})|^2$ denotes the probability density.

The MLP correction in all experiments uses 5 hidden layers, 50 neurons per layer, and \texttt{tanh} activation. For the Adam baseline, we fix the learning rate at $10^{-3}$. For the proposed NGD method, we use the adaptive step size
\begin{equation}
  \tau_k = 0.5 + 0.5\,\frac{\text{Res}(\hat{\psi})}{1 + \text{Res}(\hat{\psi})},
  \label{eq:tau_k}
\end{equation}
which accelerates early-stage descent while allowing smaller updates near convergence.
Let $R_k=\text{Res}(\hat{\psi}_k)$ denote the local variance at iteration $k$, and let $\bar R_k$ denote its exponential moving average. We choose the damping parameter based on this smoothed residual as follows:
\begin{equation}
  \lambda_k
  = \mathrm{clip}\bigl(10^{-3}\bar R_k,\;10^{-4},\;10^{-3}\bigr).
  \label{eq:lambda_k}
\end{equation}
Large residuals therefore keep $\lambda_k$ close to $10^{-3}$ during early iterations, whereas the lower bound $10^{-4}$ prevents the regularization from becoming too small near convergence. This schedule is analyzed in \Cref{sm:damping-stability}.
Unless otherwise stated, all reported main experiments use the fixed hyperparameters $(r,\alpha,K)=(128,0.8,10)$ for the Nyström rank, ADIS mixture weight, and preconditioner refresh frequency, respectively. These values are not tuned separately for individual test problems; supporting ablations are reported in \Cref{sm:nystrom-rank,sm:adis-alpha}.
Complete implementation details are available in the code repository: \url{https://github.com/cuichen1996/NGD4GPE}.

\subsection{2D optical lattice potential}
\label{subsec:exp2}

We consider the 2D optical lattice with anisotropic harmonic confinement \cite{henning2025gross}, defined by
\begin{equation}
  V(x, y)
  = \tfrac{1}{4}(x^2 + 4y^2)
  + 5\bigl(\sin^2(\pi x) + \sin^2(\pi y)\bigr).
\end{equation}
With $\beta = 250$ on $\calD = [-8,8]^2$, the ground state has a multiscale structure consisting of a macroscopic harmonic envelope modulated by high-frequency lattice oscillations.

To resolve these high-frequency features, the neural ansatz incorporates a trainable Fourier feature embedding with $\mathbf{B}\in\mathbb{R}^{16\times2}$, initialized from $\mathcal{N}(0,I_2)$, i.e., $\sigma=1.0$ in \Cref{tab:rff_config}. Expectations are estimated by M-H sampling with $4000$ samples, while the normalization constant $Z$ is computed on a uniform $100 \times 100$ grid. We compare the proposed NGD method for $100$ epochs with Adam for $2000$ epochs.

\begin{figure}[!htbp]
  \centering
  \begin{subfigure}[b]{0.32\textwidth}
    \centering
    \includegraphics[width=\textwidth]{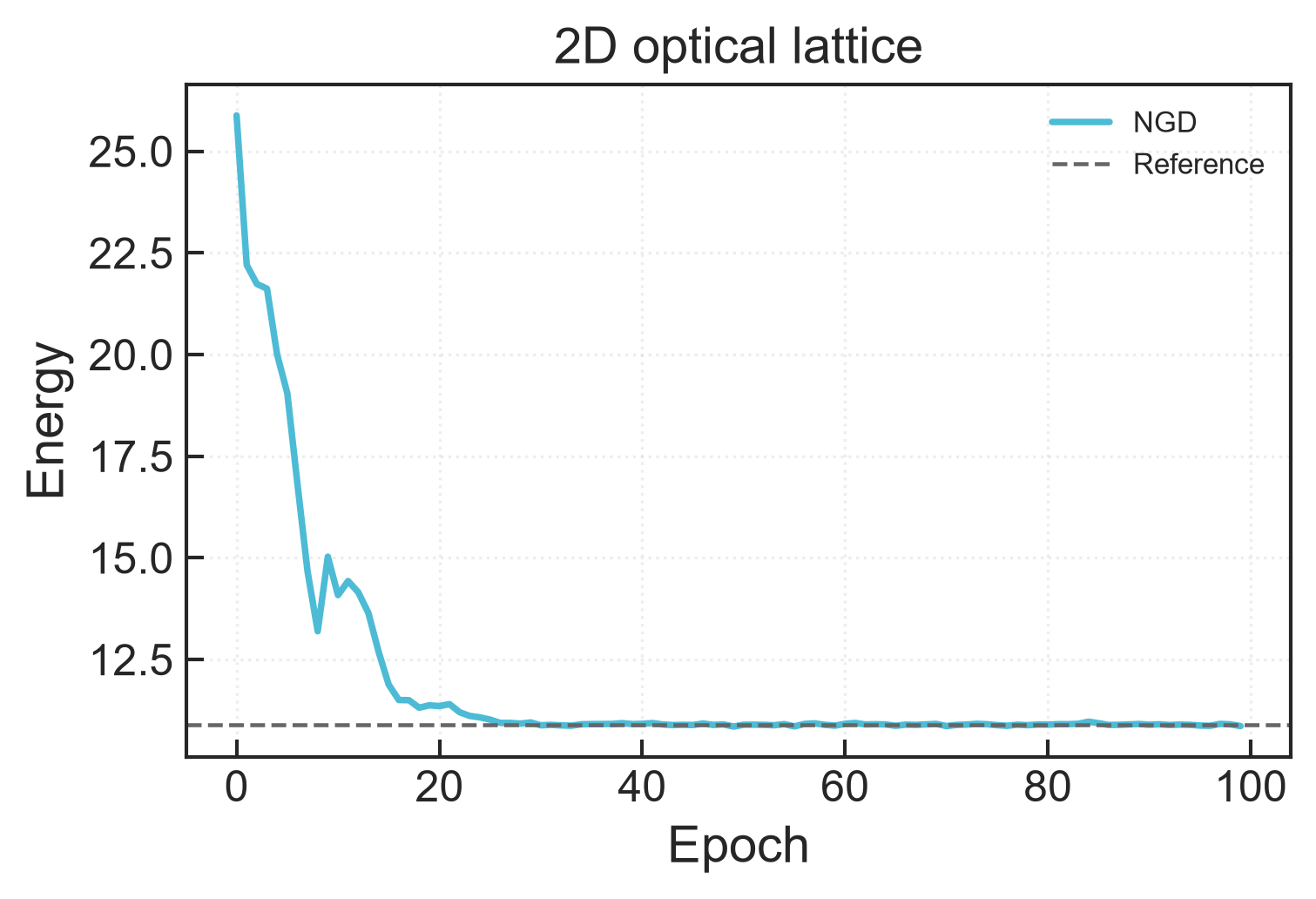}
    \caption{Evolution of the ground-state energy during NGD optimization.}
    \label{subfig:energy_lattice}
  \end{subfigure}
  \hfill
  \begin{subfigure}[b]{0.32\textwidth}
    \centering
    \includegraphics[width=\textwidth]{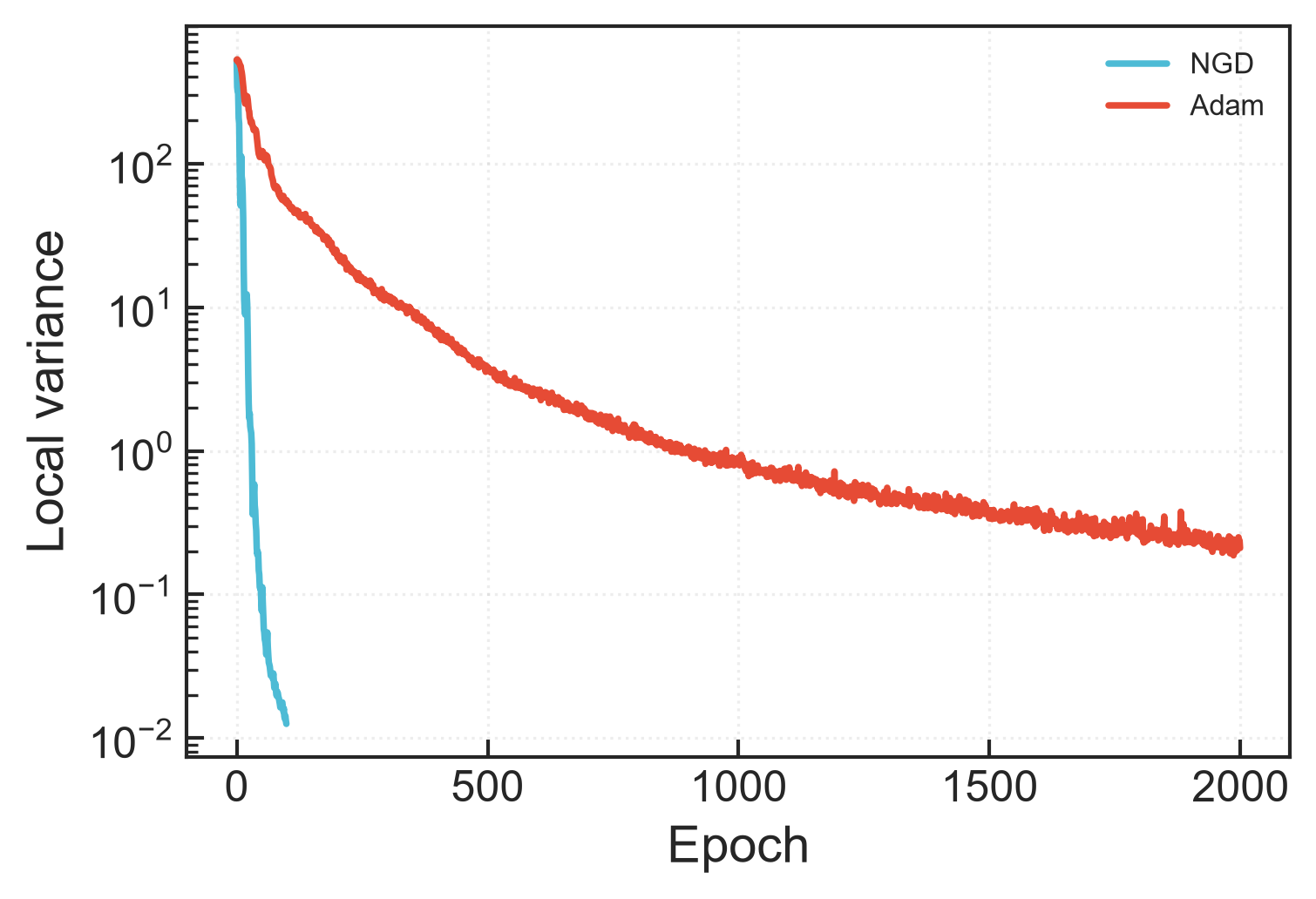}
    \caption{Comparison of Adam and NGD in terms of the local variance $\text{Res}(\hat{\psi})$.}
    \label{subfig:comparison_lattice}
  \end{subfigure}
  \hfill
  \begin{subfigure}[b]{0.32\textwidth}
    \centering
    \includegraphics[width=\textwidth]{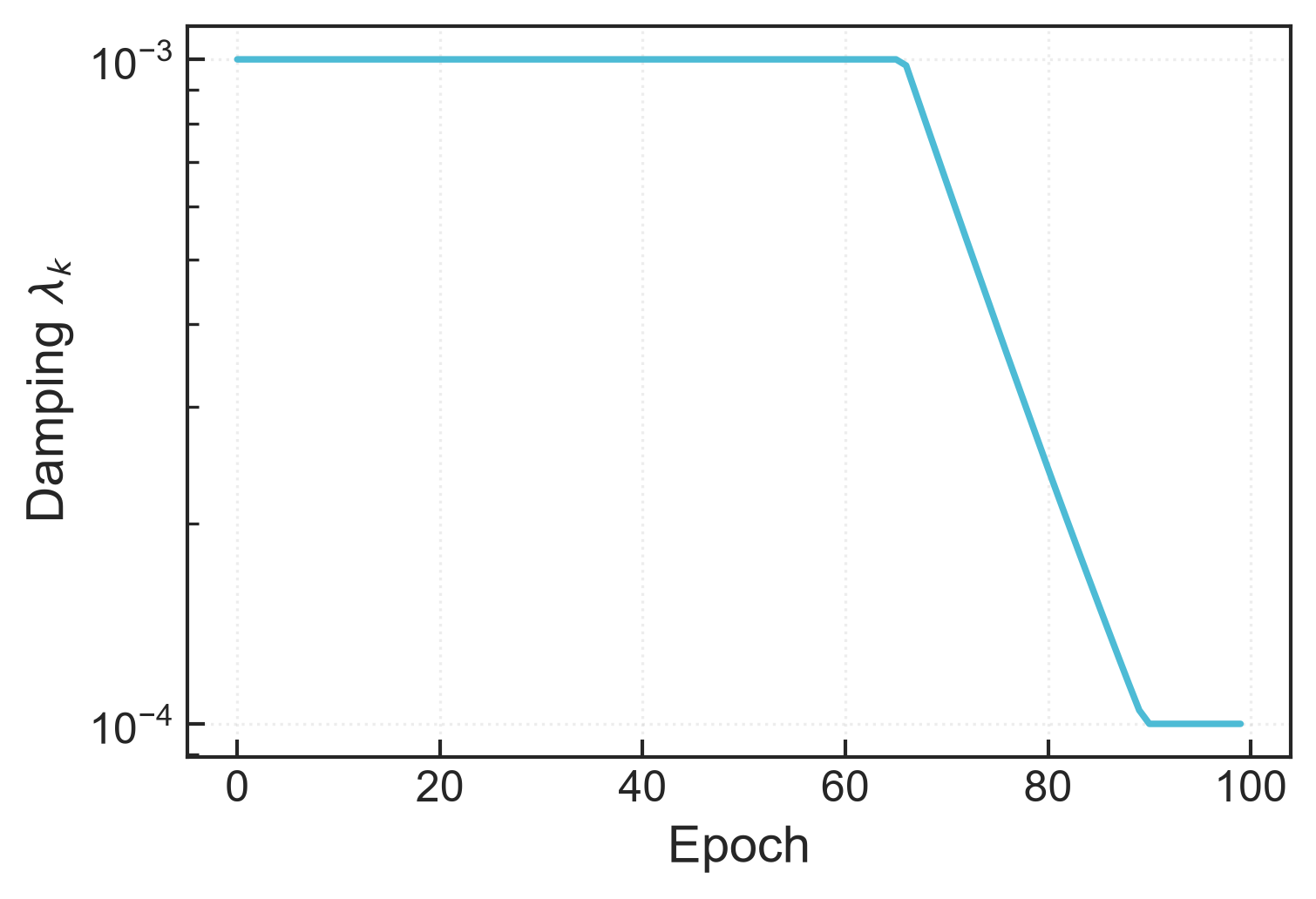}
    \caption{Adaptive damping parameter $\lambda_k$ during NGD optimization.}
    \label{subfig:lambda_lattice}
  \end{subfigure}
  \caption{Convergence history for the 2D optical lattice potential.}
  \label{fig:lattice_convergence}
\end{figure}

\Cref{fig:lattice_convergence} shows the energy evolution, local variance convergence, and adaptive damping trajectory. Adam descends rapidly at first but stagnates later. In this experiment, NGD shows stable empirical energy decay by using curvature information from the Sobolev metric, and the damping trajectory illustrates how the residual-aware schedule regularizes the noisy early phase before gradually releasing curvature information.

\begin{figure}[!htbp]
  \centering
  \includegraphics[width=\textwidth]{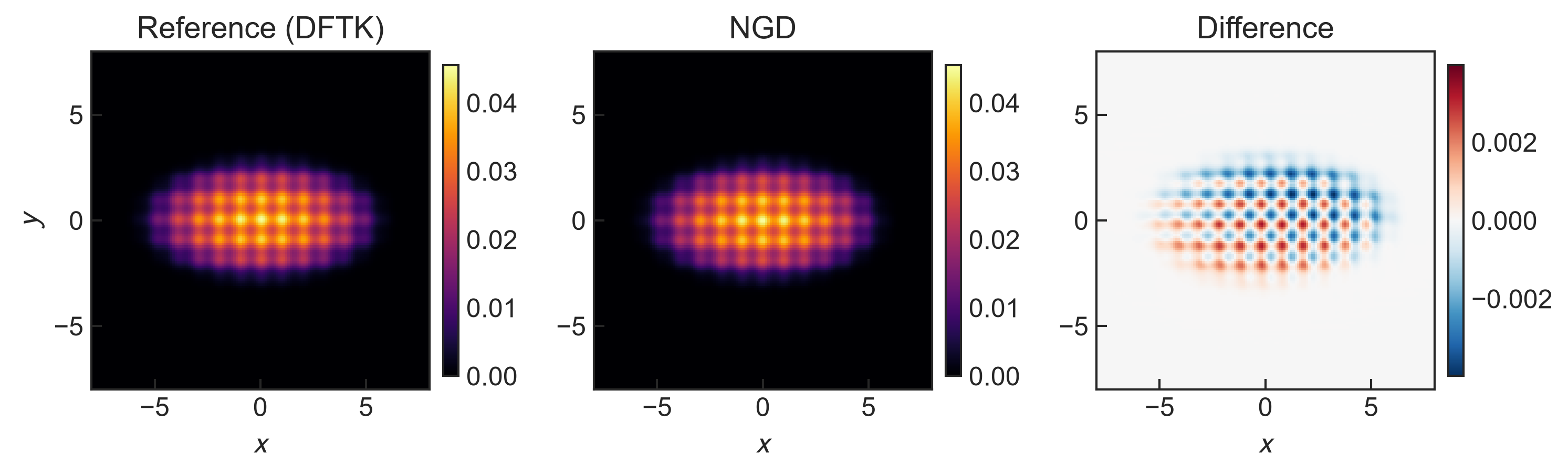}
  \caption{Ground-state density for the 2D optical lattice potential computed by NGD.}
  \label{fig:lattice_results}
\end{figure}

\Cref{fig:lattice_results} visualizes the ground-state density obtained by NGD.
The solution captures both the smooth global envelope and the checkerboard-like periodic modulation induced by the lattice potential, showing that the method resolves features across multiple spatial scales.

\Cref{tab:lattice_results} reports the quantitative error metrics at convergence.
Both optimizers reach comparable accuracy, but NGD needs far fewer iterations.
\begin{table}[!htbp]
  \centering
  \caption{Error comparison for the 2D optical lattice problem.}
  \label{tab:lattice_results}
    \begin{tabular}{ccccc}
      \toprule
      Optimizer & $\varepsilon_E$ & $\varepsilon_\mu$ & $\varepsilon_\rho$ & Iter. \\
      \midrule
      Adam & $2.90 \times 10^{-3}$ & $5.89 \times 10^{-4}$ & $6.31 \times 10^{-2}$ & 1986 \\
      NGD  & $3.46 \times 10^{-3}$ & $1.39 \times 10^{-4}$ & $6.22 \times 10^{-2}$ & \textbf{99} \\
      \bottomrule
  \end{tabular}
\end{table}

Finally, \Cref{tab:resource_comparison} compares the computational cost of different optimization and preconditioning strategies.
NGD (Explicit $G$) uses the full Gram matrix without preconditioning. NGD (Block-Jacobi-PCG) uses a layer-wise block-Jacobi preconditioner, which acts as a nonoverlapping domain-decomposition preconditioner in parameter space. NGD (Nyström-PCG) is the proposed Nyström-preconditioned approach. Although NGD costs more per iteration than Adam, faster convergence reduces total training time. Nyström-PCG reduces GPU memory by more than an order of magnitude compared with the explicit Gram matrix and achieves a $3.6\times$ end-to-end speedup over Adam.

\begin{table}[!htbp]
  \centering
  \caption{Comparison of computational costs for different optimization and preconditioning strategies on an A100 GPU.}
  \label{tab:resource_comparison}
  \resizebox{\textwidth}{!}{
    \begin{tabular}{lccccc}
      \toprule
      Algorithm & $\varepsilon_\rho$ & Epochs & Time (s) & Speedup & VRAM (GB) \\
      \midrule
      Adam (Baseline) & $7.2\%$ & 2000 & 67.1 & 1.0$\times$ & 0.7 \\
      NGD (Explicit $G$) & $6.2\%$ & 100 & 61.9 & 1.1$\times$ & 25.1 \\
      NGD (Block-Jacobi-PCG) & $6.3\%$ & 100 & 23.5 & 2.8$\times$ & 4.1 \\
      \textbf{NGD (Nyström-PCG)} & $\mathbf{6.3\%}$ & 100 & \textbf{18.7} & \textbf{3.6$\times$} & 2.4 \\
      \bottomrule
  \end{tabular}}
\end{table}

\subsubsection{\texorpdfstring{Comparison with the $L^{2}$ natural gradient (SR)}{Comparison with the L2 natural gradient (SR)}}
\label{subsec:sr_comparison}

\begin{figure}[!htbp]
  \centering
  \begin{subfigure}[t]{0.27\textwidth}
    \centering
    \includegraphics[width=\textwidth]{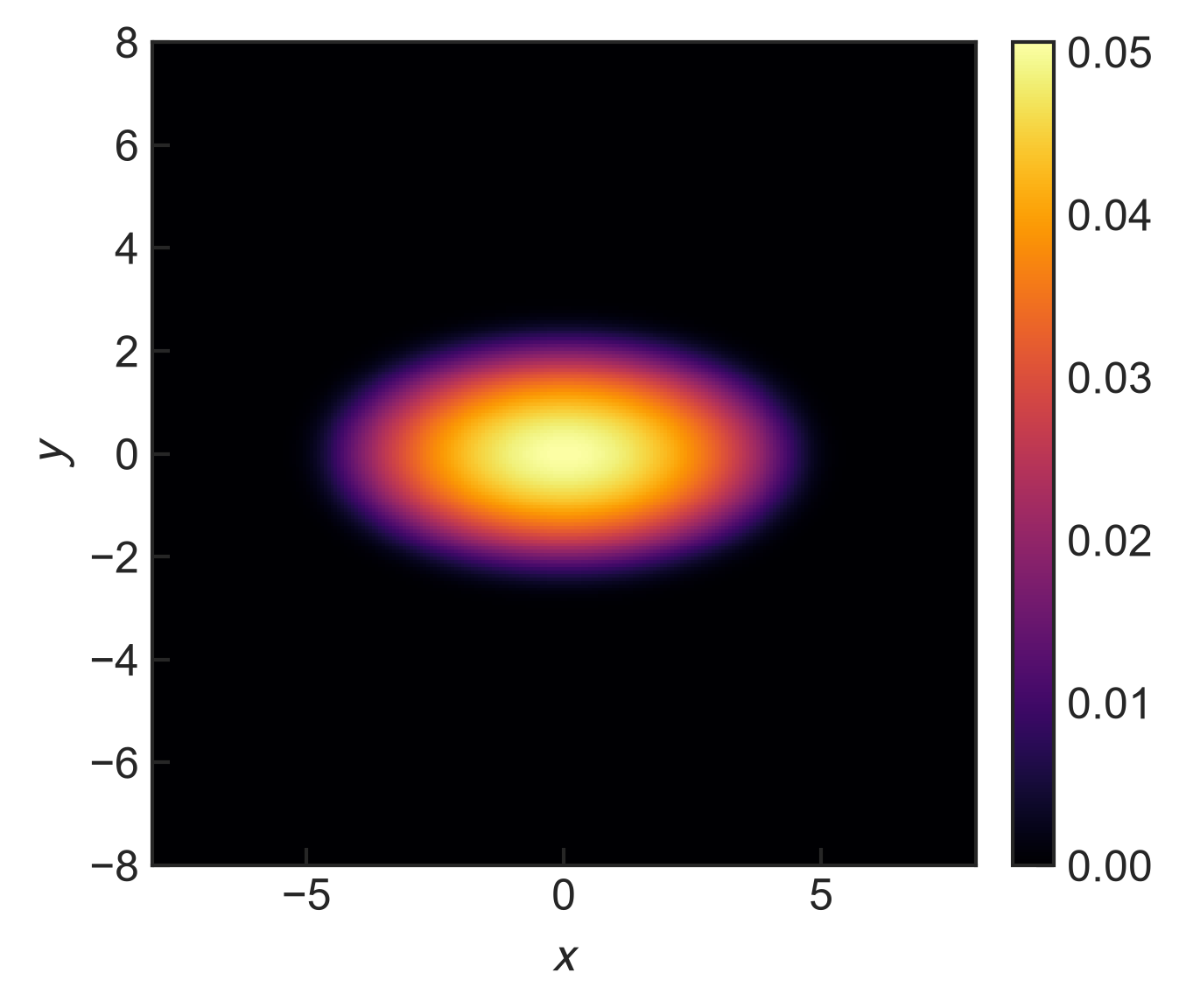}
    \caption{Harmonic density.}
  \end{subfigure}\hfill
  \begin{subfigure}[t]{0.34\textwidth}
    \centering
    \includegraphics[width=\textwidth]{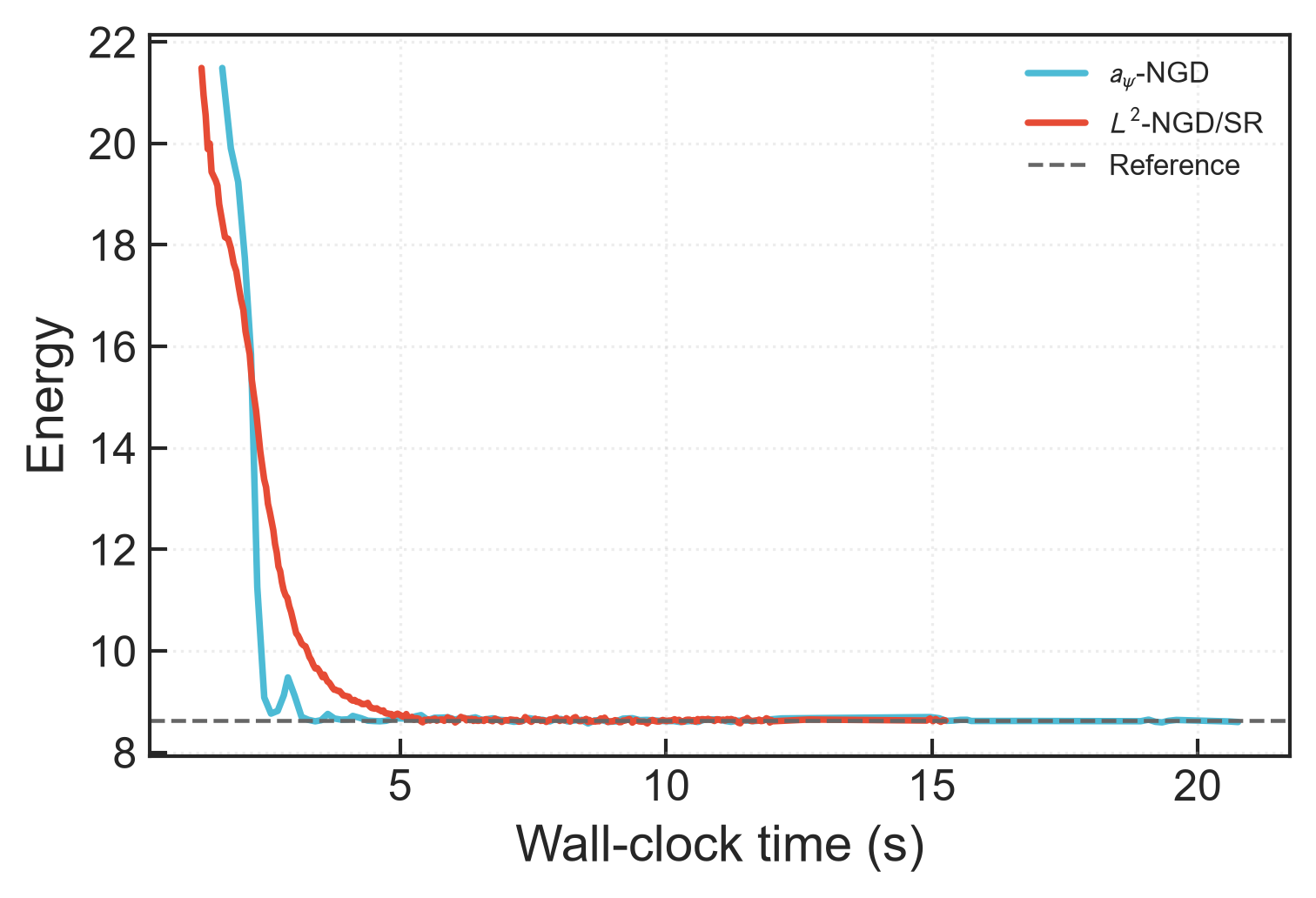}
    \caption{Harmonic energy.}
  \end{subfigure}\hfill
  \begin{subfigure}[t]{0.34\textwidth}
    \centering
    \includegraphics[width=\textwidth]{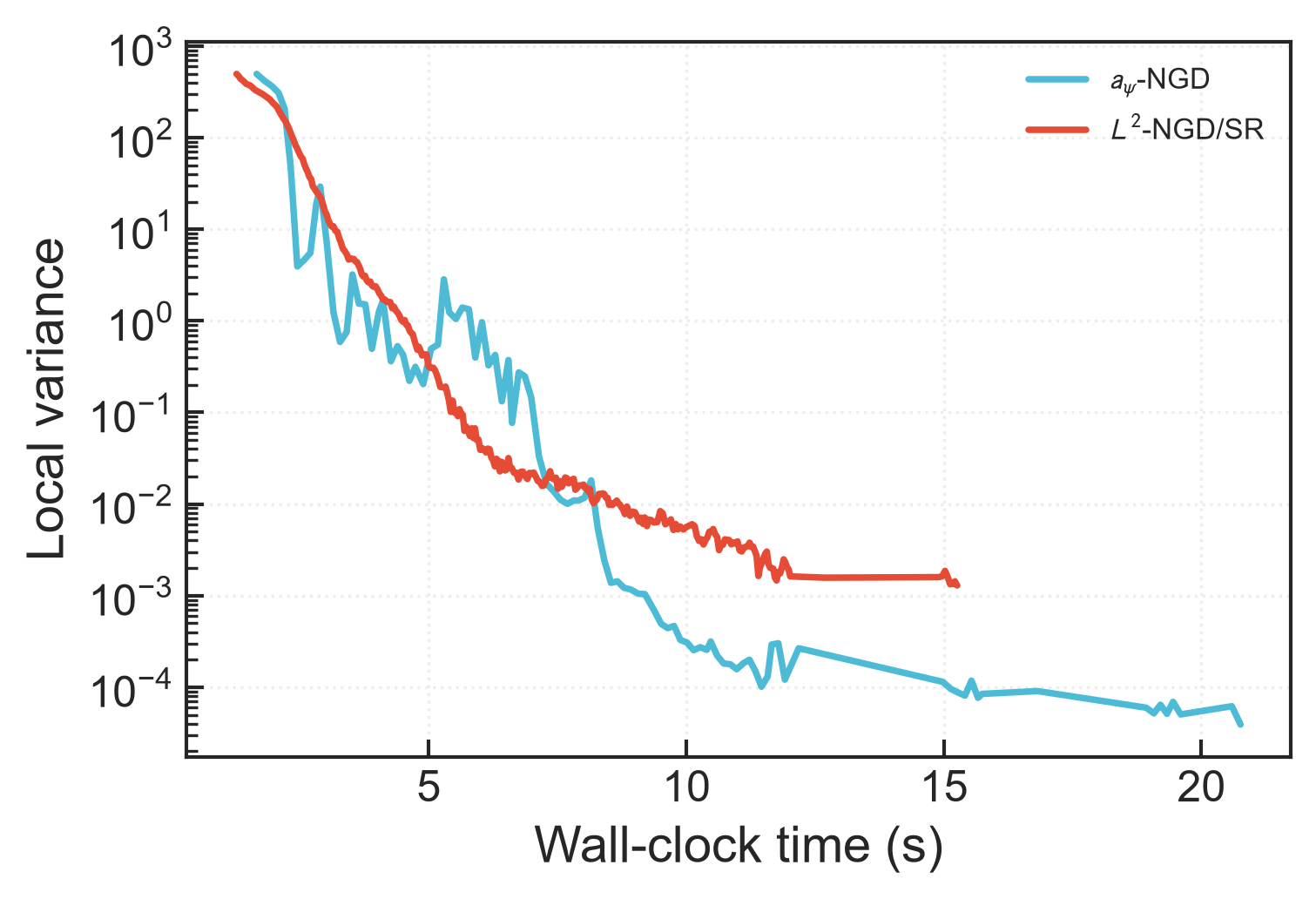}
    \caption{Harmonic variance.}
  \end{subfigure}

  \begin{subfigure}[t]{0.27\textwidth}
    \centering
    \includegraphics[width=\textwidth]{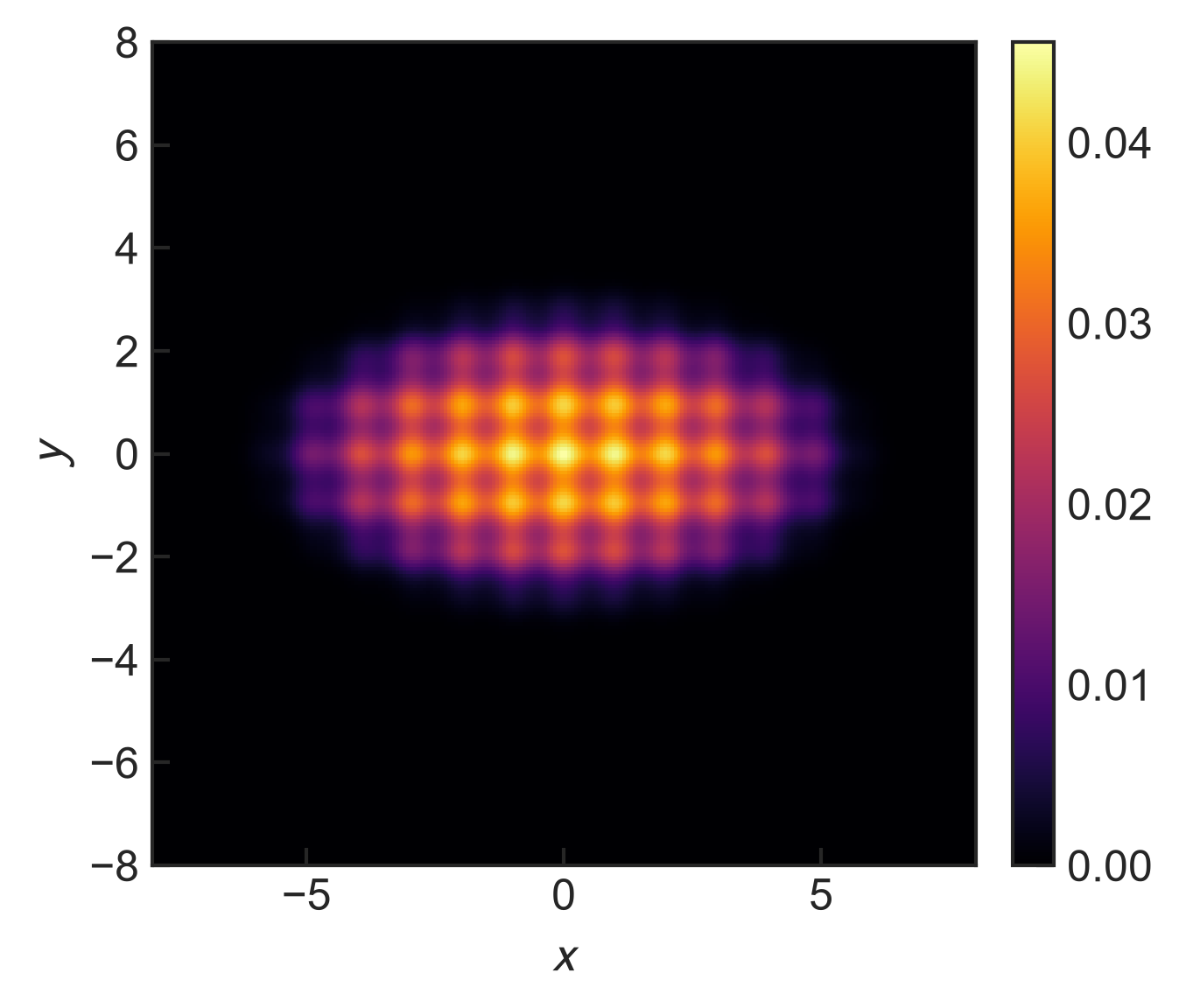}
    \caption{Lattice density.}
  \end{subfigure}\hfill
  \begin{subfigure}[t]{0.34\textwidth}
    \centering
    \includegraphics[width=\textwidth]{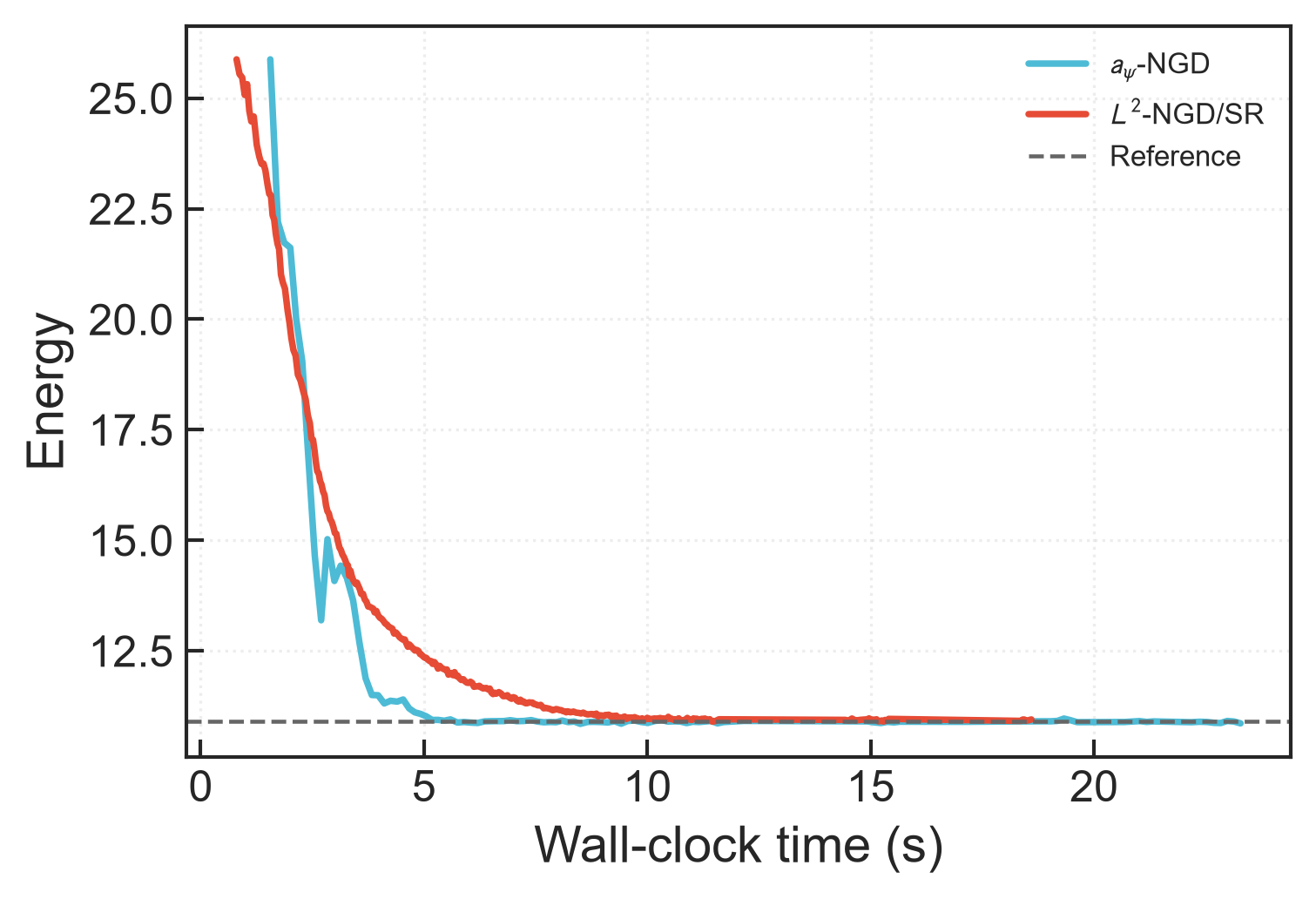}
    \caption{Lattice energy.}
  \end{subfigure}\hfill
  \begin{subfigure}[t]{0.34\textwidth}
    \centering
    \includegraphics[width=\textwidth]{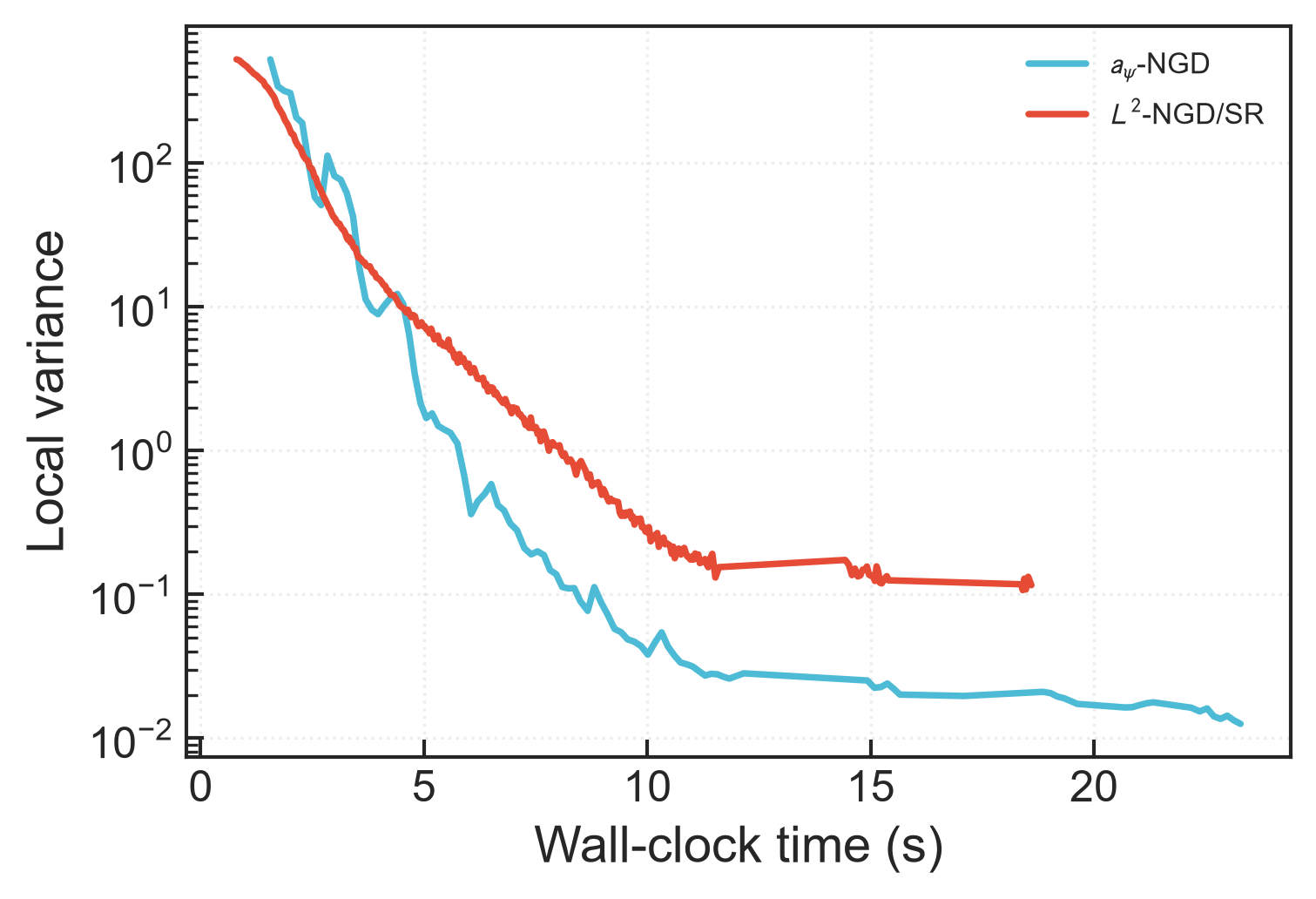}
    \caption{Lattice variance.}
  \end{subfigure}
  \caption{Comparison of $a_\psi$-NGD and $L^{2}$-NGD/SR on two nonlinear 2D GPE problems with $\beta=250$. The first and second rows show the harmonic trap and optical lattice, respectively; columns show the ground-state density, energy versus wall-clock time, and local variance versus wall-clock time.}
  \label{fig:sr_comparison}
\end{figure}

We next report a wall-clock comparison on two nonlinear 2D ground-state problems: a harmonic trap with $\beta=250$ and an optical lattice with $\beta=250$. Within each problem, both methods use $N=4000$ MCMC samples and the same five-hidden-layer network of width $50$. The $a_\psi$-NGD method is run for $100$ epochs and the $L^{2}$-NGD/SR method for $300$ epochs, all on a single A100 GPU. \Cref{fig:sr_comparison} reports the ground-state density, energy, and local variance as functions of wall-clock time. In both tests, $a_\psi$-NGD reaches a substantially lower local variance despite using one third as many epochs, demonstrating an end-to-end wall-clock advantage rather than only better per-step progress.

\subsection{3D optical lattice potential}
\label{subsec:3D_benchmark}

To validate the ADIS strategy, we use a 3D optical lattice benchmark for which deterministic grid-based quadrature is still feasible and provides a reliable reference. Following \cite{bao2025computing}, the external trapping potential is
\begin{equation}
  V(x,y,z) = \frac{1}{2}(x^2+y^2+z^2) + \frac{5}{2} \left(\sin^2\!\left(\frac{\pi x}{4}\right) + \sin^2\!\left(\frac{\pi y}{4}\right) + \sin^2\!\left(\frac{\pi z}{4}\right)\right),
  \label{eq:3d_pot}
\end{equation}
with interaction strength $\beta = 200$ and computational domain $\calD = [-6,6]^3$.
We compare two strategies for evaluating the normalization constant $Z$: a deterministic Riemann sum on a uniform tensor grid with $100^3$ points, and the proposed ADIS method using $8000$ samples with a mixture ratio $\alpha = 0.8$.
Because the potential contains the lattice frequency $\pi/4$ in each coordinate, this experiment uses a trainable Fourier feature embedding with $\mathbf{B}\in\mathbb{R}^{16\times3}$ initialized from $\mathcal{N}(0,0.25^2I_3)$, as summarized in \Cref{tab:rff_config}.

\begin{figure}[!htbp]
  \centering
  \begin{subfigure}[b]{0.36\textwidth}
    \centering
    \includegraphics[width=\textwidth]{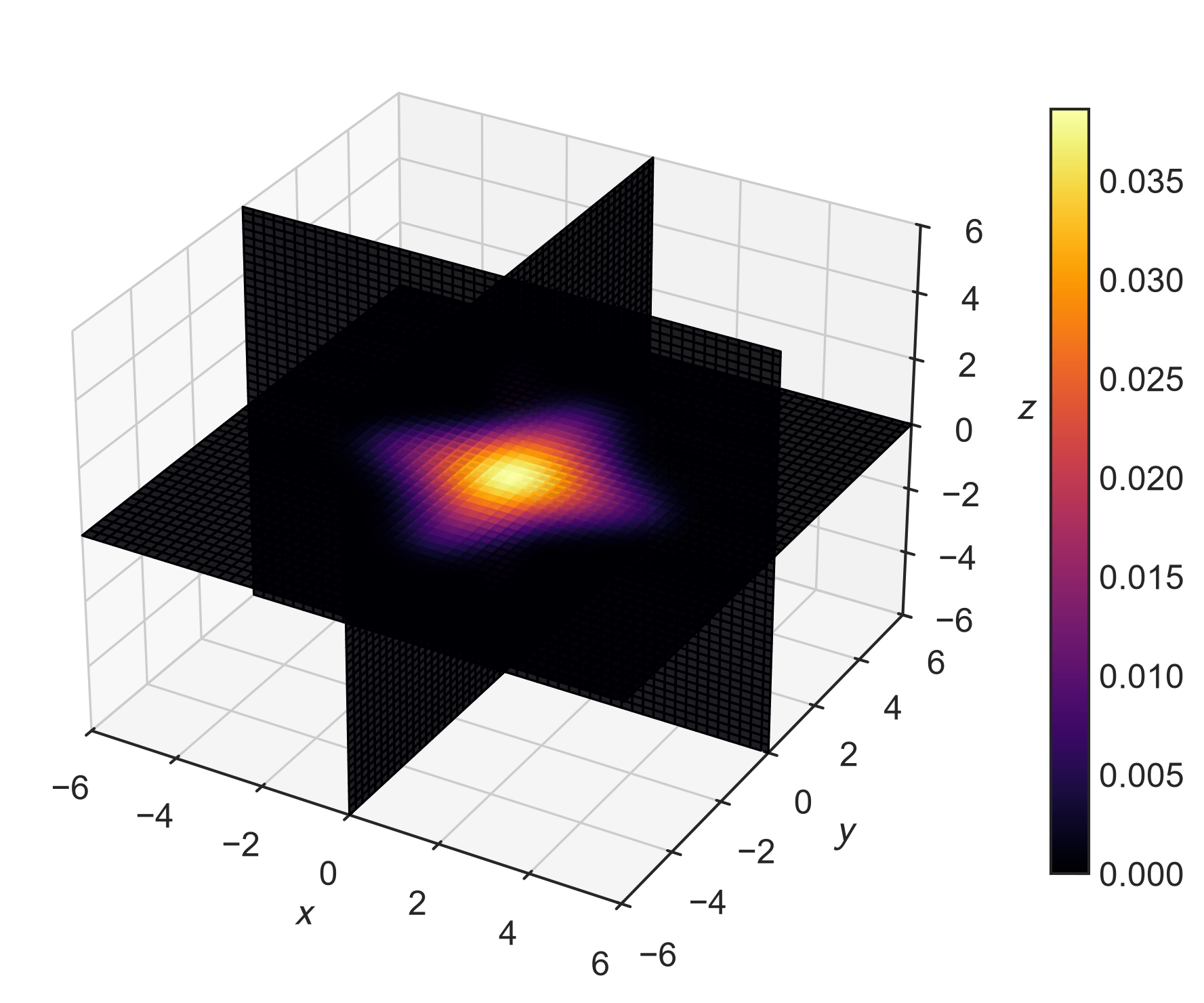}
    \caption{Planar density slices.}
    \label{fig:density3D}
  \end{subfigure}
  \hfill
  \begin{subfigure}[b]{0.30\textwidth}
    \centering
    \includegraphics[width=\textwidth]{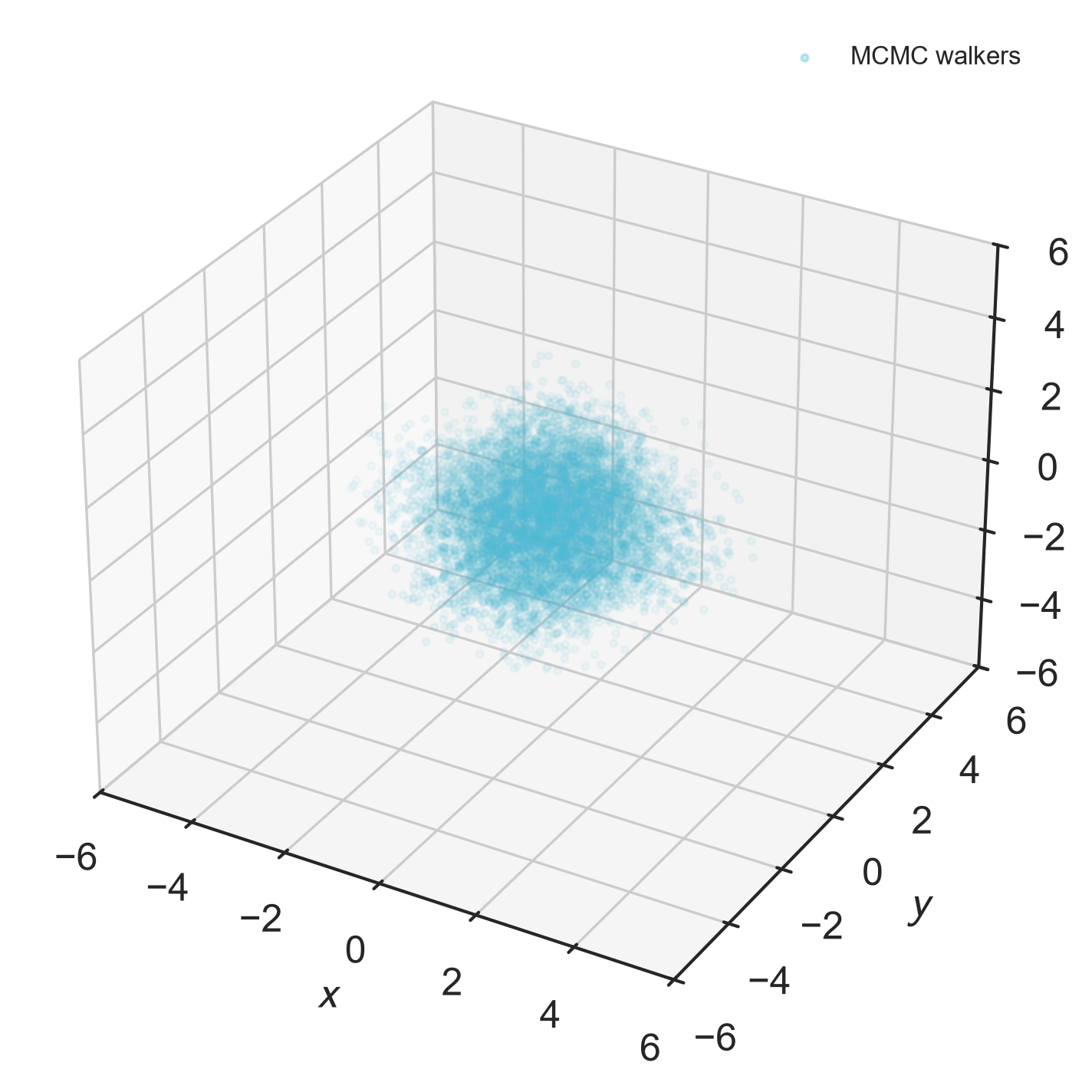}
    \caption{MCMC walker distribution.}
    \label{fig:mcmc}
  \end{subfigure}
  \hfill
  \begin{subfigure}[b]{0.30\textwidth}
    \centering
    \includegraphics[width=\textwidth]{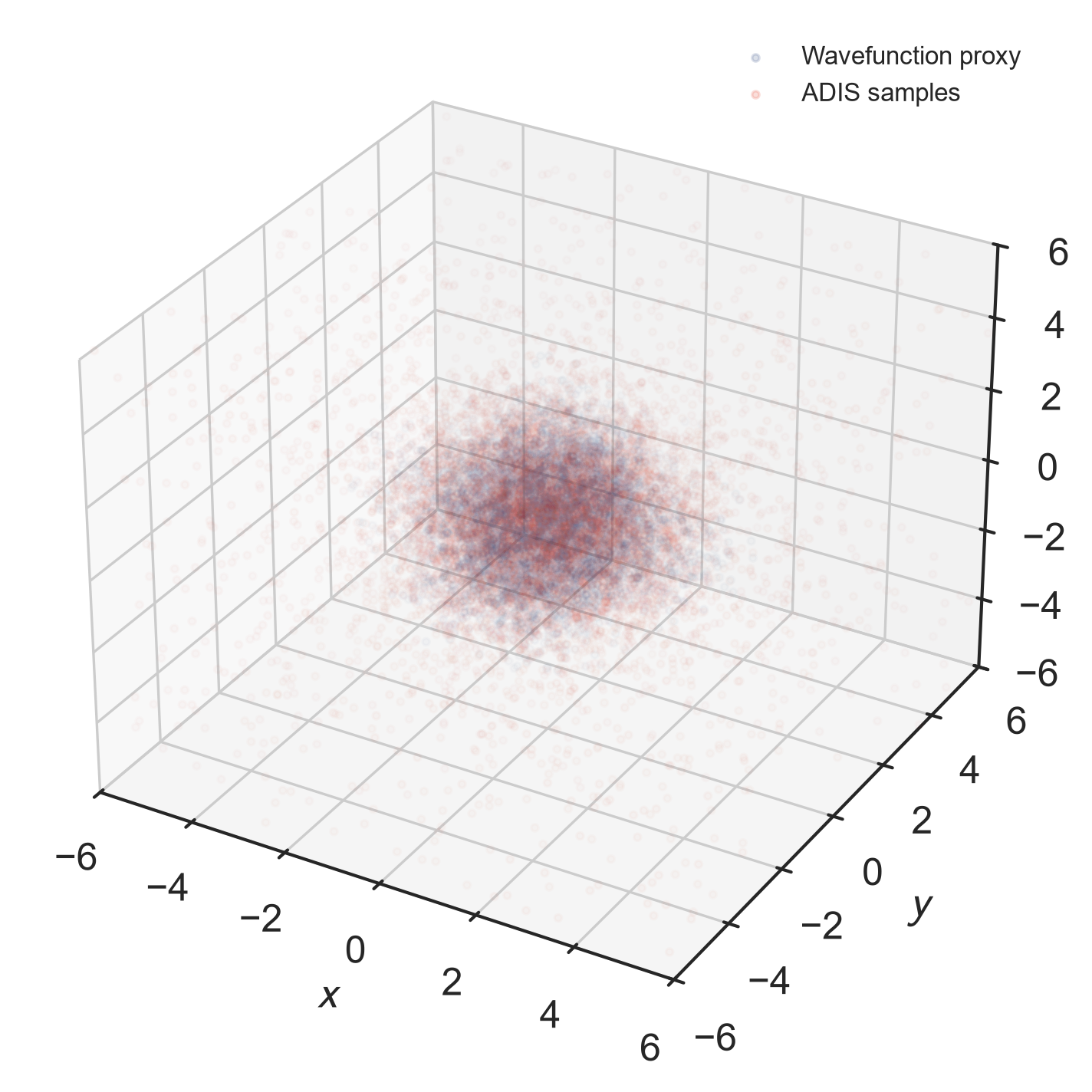}
    \caption{ADIS integration samples.}
    \label{fig:z_sample}
  \end{subfigure}
  \caption{Visualization of the ground-state density and sampling profiles for the 3D optical lattice problem.}
  \label{fig:model1_visual}
\end{figure}

\Cref{fig:model1_visual} shows both the neural VMC solution and the sampling behavior induced by ADIS.
Specifically, \Cref{fig:density3D} displays planar slices of the  ground-state density, \Cref{fig:mcmc} shows the spatial distribution of MCMC walkers during training, and \Cref{fig:z_sample} shows the ADIS samples used to estimate $Z$.
The ADIS sample distribution closely matches the ground-state density, indicating that the adaptive proposal concentrates samples in high-probability regions.

\Cref{tab:3D_adis} compares the two computational strategies.
ADIS achieves accuracy comparable to that of the dense grid-based Riemann sum while using only about $1\%$ of the sample points. Because grid-based quadrature costs grow exponentially with dimension, ADIS avoids the exponential sample and storage cost of full tensor-product quadrature. Under this experimental configuration, the two strategies have comparable peak allocated memory.
\begin{table}[!htbp]
  \centering
  \caption{Comparison of accuracy and computational cost between deterministic grid-based integration and the ADIS framework for the 3D optical lattice problem.}
  \label{tab:3D_adis}
    \begin{tabular}{lcccc}
      \toprule
      Method & Sample points & Time (s) & $\varepsilon_\rho$ & VRAM (GB) \\
      \midrule
      Riemann sum & $100^3$ & 42.6 & $5.82 \times 10^{-2}$ & 5.986 \\
      \textbf{ADIS} & $8000$ & 45.0 & $5.82 \times 10^{-2}$ & 5.964 \\
      \bottomrule
  \end{tabular}
\end{table}

\subsubsection{Neural-network solution as an initial guess}

We also use the NGD solution as an initial guess for a high-precision solver. The neural wavefunction is evaluated on a uniform real-space grid and then projected onto a plane-wave basis by FFT. This spectral representation initializes the energy-adaptive Riemannian conjugate gradient (EARCG) method \cite{peterseim2025energy}. In practice, 20 NGD iterations, requiring less than 10 seconds, suffice as a warm start for refinement.

\begin{figure}[!htbp]
  \centering
  \includegraphics[width=0.8\textwidth]{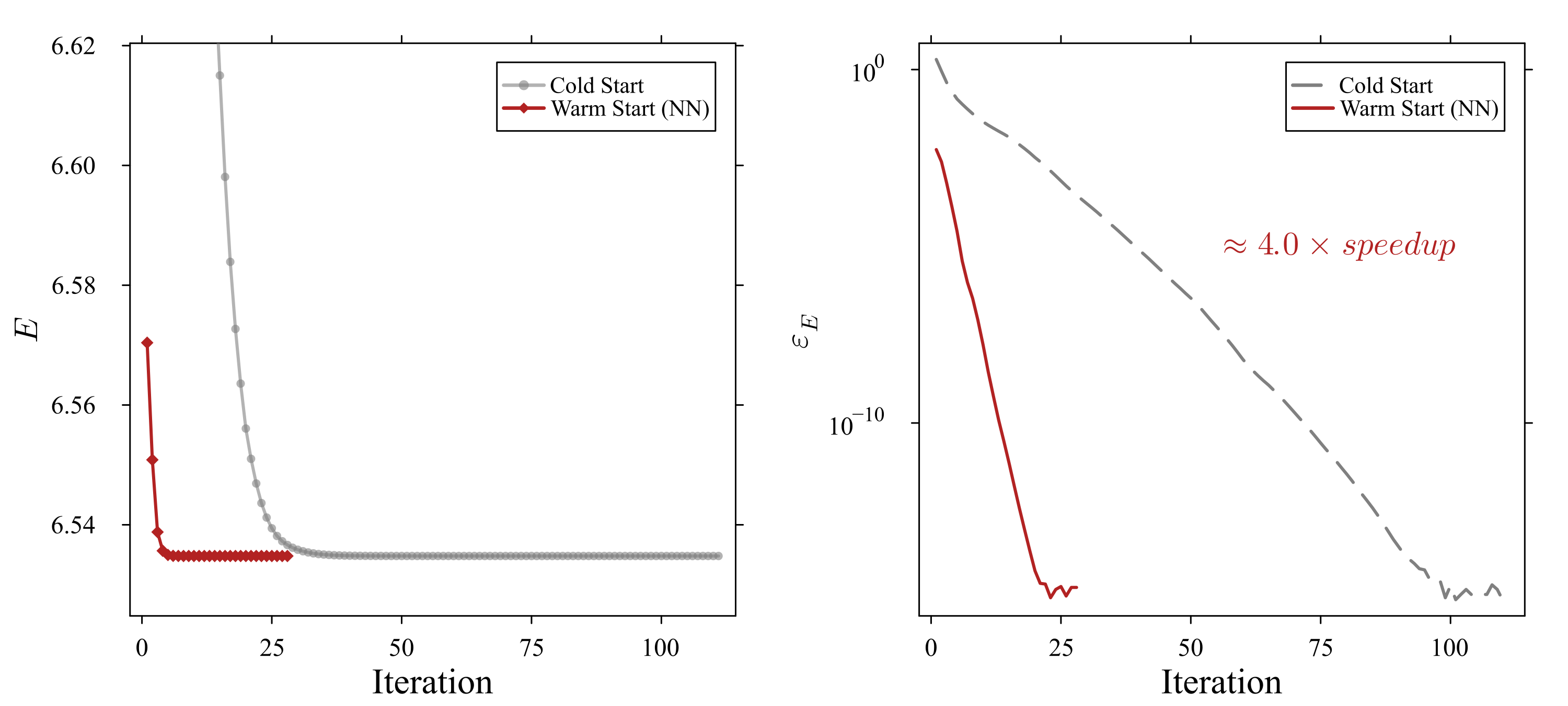}
  \caption{Energy convergence of the EARCG solver with a neural-network warm start versus a default initial guess.}
  \label{fig:initialization}
\end{figure}
As \Cref{fig:initialization} shows, this warm start improves over a default initialization. It skips the initial high-energy relaxation phase and reduces the starting energy error. The EARCG iterations then converge faster on a logarithmic scale, indicating that the neural network has already captured the global structure of the ground state. The warm-started iterate is therefore closer to the convergence basin, reducing the ill-conditioning typical of nonlinear GPE landscapes.

\subsection{Comparison with the norm-DNN}
\label{subsec:comparison_normdnn}

We now compare with the norm-DNN approach \cite{bao2025computing}, a PINN-based method that enforces the mass constraint through an explicit normalization layer.
In the 3D setting, norm-DNN considers the discretized optimization problem
\begin{equation}
  \min_{\theta}
  \left\{
    \mathcal{L}
    =
    \frac{|\calD|}{N_{b}}
    \sum_{j,k,l}
    \left[
      \frac{1}{2}\left|\nabla \hat{\psi}_\theta\right|^2
      + V\left|\hat{\psi}_\theta\right|^2
      + \frac{1}{2}\beta\left|\hat{\psi}_\theta\right|^4
    \right]_{(x_j,y_k,z_l)}
  \right\}.
\end{equation}

We compare the methods on two benchmark problems from \cite{bao2025computing}, following the original experimental settings.
The first example is a 2D optical lattice problem with potential
\begin{equation}
  V(x,y)
  =
  \frac{1}{2}\left(x^2+y^2\right)
  +
  \frac{5}{2}
  \left(
    \sin^2\!\left(\frac{\pi x}{4}\right)
    +
    \sin^2\!\left(\frac{\pi y}{4}\right)
  \right),
\end{equation}
posed on the domain $\calD=[-6,6]^2$ with interaction strength $\beta=400$.
The second example is the 3D optical lattice problem defined by \eqref{eq:3d_pot} on $\calD=[-6,6]^3$ with $\beta=200$.

\Cref{tab:normdnn_comparison} summarizes the results. Our method converges in 100 epochs, while norm-DNN requires 30,000 to 100,000 epochs for comparable accuracy. For the 3D optical lattice, our method achieves a lower $L^2$ density error with far fewer sampling points.
\begin{table}[!htbp]
  \centering
  \caption{Comparison between the projected Sobolev NGD method and the norm-DNN method \cite{bao2025computing} for 2D and 3D optical lattice benchmark problems.}
  \label{tab:normdnn_comparison}
  \renewcommand{\arraystretch}{1.35}
  \setlength{\tabcolsep}{6pt}
  \scriptsize
  \resizebox{\textwidth}{!}{
    \begin{tabular}{lccccc}
      \toprule
      Method &
      Batch size &
      Epochs &
      $\varepsilon_\rho$ &
      $\varepsilon_E$ &
      Time (s) \\
      \midrule
      \multicolumn{6}{c}{\textbf{2D case}} \\
      \midrule
      norm-DNN & 4096 & 30{,}000 & $4.82\times10^{-2}$ & $2.08\times10^{-3}$ & 169.7 \\
      \textbf{Ours} & 4000 & \textbf{100} & $\mathbf{3.14\times10^{-2}}$ & $\mathbf{7.51\times10^{-4}}$ & \textbf{15.1} \\
      \midrule
      \multicolumn{6}{c}{\textbf{3D case}} \\
      \midrule
      norm-DNN & 60{,}000 & 100{,}000 & $1.28\times10^{-1}$ & $1.60\times10^{-2}$ & 687.9 \\
      \textbf{Ours} & 8000 & \textbf{100} & $\mathbf{5.94\times10^{-2}}$ & $\mathbf{4.92\times10^{-3}}$ & \textbf{39.4} \\
      \bottomrule
\end{tabular}}
\end{table}

\subsection{High-dimensional scalability test}
\label{subsec:highd_scaling}

To assess scalability, we use a fully coupled high-dimensional nonlinear Schrödinger model with an anisotropic quadratic potential,
\begin{equation}
  V(\mathbf{r}) = \frac{1}{2}\,\mathbf{r}^\top \mathbf{A}\mathbf{r},
\end{equation}
where $\mathbf{A} \in \mathbb{R}^{d\times d}$ is a randomly generated symmetric positive definite matrix with condition number fixed at $\kappa(\mathbf A)=10$ in all dimensions.
The full matrix $\mathbf A$ couples all spatial dimensions and precludes a separation-of-variables approach. We set $\beta = 2000$ on $\calD = [-6,6]^d$. Grid-based methods are infeasible here because the number of degrees of freedom grows exponentially with $d$.

Because high-dimensional reference solutions are unavailable, we use the generalized virial theorem \cite{dalfovo1999theory} as a self-consistency check. We define the virial residual ratio
\begin{equation}
  \mathcal{R}_{\mathrm{vir}} = \frac{2\langle T \rangle + d\langle E_{\mathrm{int}} \rangle}{2\langle V \rangle},
\end{equation}
where $\langle T \rangle$, $\langle V \rangle$, and $\langle E_{\mathrm{int}} \rangle$ denote the kinetic, potential, and interaction energy components, respectively. For an exact normalized ground state, this ratio equals one. Deviations from one therefore measure the physical inconsistency of the numerical approximation.

We consider dimensions $d=4,\ldots,8$. Each model is trained for $200$ epochs using $15{,}000$ HMC walkers, while $12{,}000$ ADIS samples are used to estimate the normalization constant $Z$.
The MLP correction has five hidden layers of width $50$, with both depth and width held fixed across dimensions. Raw coordinates are passed directly to the network without an RFF embedding, consistent with the smooth, low-frequency structure of this benchmark. For the Gaussian envelope, the implementation registers a full $d\times d$ matrix $\mathbf C$, distinct from the potential matrix $\mathbf A$, but uses only its lower-triangular part in the forward evaluation. Counting all MLP weights and biases together with all $d^2$ registered entries of $\mathbf C$ gives a total registered parameter count of $P=d^2+50d+10{,}301$. Consequently, the input layer and the registered envelope matrix are the only dimension-dependent components, so the increase in $P$ remains modest over the tested range.
All high-dimensional runs evaluate the Sobolev Gram-matrix products using the implicit operator described in \Cref{subsubsec:matrix-free-implementation}; consequently, neither the centered score matrix $\mathbf J$ nor the mixed input--parameter derivative tensor $\mathbf K$ is materialized.
The corresponding parameter counts and numerical results are summarized in \Cref{tab:hd_verify}.
Across the tested dimensions, the virial ratio remains close to its exact value and the local variance stays small. The measured wall-clock time varies only modestly, and the peak allocated GPU memory is nearly constant. Thus, with fixed MLP width and depth and only mild parameter-count growth, the observed computational cost shows no evident exponential growth over this benchmark.
\begin{table}[!htbp]
  \centering
  \caption{Architecture size,  physical consistency, and computational cost for the high-dimensional fully coupled model.}
  \label{tab:hd_verify}
  \resizebox{\textwidth}{!}{
    \begin{tabular}{cccccc}
      \toprule
      Dimension & $P$ & Time (s) & Peak allocated VRAM (GB) & Local variance & Virial ratio \\
      \midrule
      4 & 10,517 & 174.3 & 0.255 & $1.075\times10^{-4}$ & 1.00267 \\
      5 & 10,576 & 170.2 & 0.256 & $2.556\times10^{-4}$ & 0.99735 \\
      6 & 10,637 & 158.5 & 0.258 & $2.993\times10^{-4}$ & 0.99998 \\
      7 & 10,700 & 172.7 & 0.260 & $1.728\times10^{-4}$ & 1.00293 \\
      8 & 10,765 & 162.4 & 0.261 & $1.251\times10^{-4}$ & 1.00681 \\
      \bottomrule
  \end{tabular}}
\end{table}

\section{Conclusion and outlook}
\label{sec:conclusion}

This paper presents a projected Sobolev NGD method for computing high-dimensional GPE ground states. The framework maps the continuous gradient flow onto the neural tangent space by Galerkin projection, preserving normalization and geometric consistency. To handle the GPE's nonlinearity and dimensionality, we developed a hybrid ADIS--MCMC sampling strategy and a matrix-free Nyström-preconditioned PCG solver for Gram matrix inversion. The method converges an order of magnitude faster than PINN approaches and the implicit Sobolev Gram-operator implementation avoids both the $\mathcal{O}(NP)$ score-matrix cache and the $\mathcal{O}(NPD)$ mixed-derivative tensor, with nearly constant measured allocated memory over the tested $d=4,\ldots,8$ benchmark. The learned solutions also provide warm starts that accelerate traditional solvers.

Several directions remain for future work. Extending the method to rotating GPEs will require complex-valued neural networks capable of representing vortices and metastable states. Computing excited states using penalty or deflation strategies is another natural direction. For dimensions in which a fixed hidden width is insufficient, a dimension-dependent width rule $P=P(d)$ will be needed. The present fixed-width test should therefore not be interpreted as a general expressivity result. Spectral bias in high-frequency regimes, such as Anderson localization, also motivates architectures beyond standard MLPs. For high-dimensional multiscale or quasi-periodic problems, i.i.d. Gaussian RFFs may be inefficient. Structured frequency dictionaries, separable or tensorized Fourier features, and multiscale wavelet-type encodings are promising alternatives.

\appendix
\section{Stability interpretation of adaptive damping}
\label[appendix]{sm:damping-stability}

This appendix records the conditional-expectation estimate underlying the adaptive damping discussion in the main text. The estimate provides a stability rationale for the regularized stochastic update, but it is not a pathwise monotonicity theorem for the realized iterates.

Let
\[
\mathbf{g}_k=\nabla_{\boldsymbol{\theta}} L(\boldsymbol{\theta}^k),
\qquad
\mathbf{G}_k=\mathbf{G}(\boldsymbol{\theta}^k),
\qquad
\mathbf{H}_k=\mathbf{G}_k+\lambda_k\mathbf{I},
\]
where $L(\boldsymbol{\theta})=E(\hat{\psi}_{\boldsymbol{\theta}})$ is the parameter-space loss and $\lambda_k>0$ is the damping parameter. We write the Monte Carlo gradient estimator as
\[
\widehat{\mathbf{g}}_k=\mathbf{g}_k+\boldsymbol{\xi}_k,
\qquad
\mathbb{E}[\boldsymbol{\xi}_k\mid\mathcal{F}_k]=0,
\qquad
\mathbb{E}[\boldsymbol{\xi}_k\boldsymbol{\xi}_k^\top\mid\mathcal{F}_k]=\boldsymbol{\Sigma}_k,
\qquad
\operatorname{tr}(\boldsymbol{\Sigma}_k)\le\sigma_k^2 .
\]
Here $\mathcal{F}_k$ is the information available before Monte Carlo sampling at iteration~$k$, $\boldsymbol{\Sigma}_k$ is the conditional covariance of the gradient noise, and $\sigma_k^2$ is a scalar upper bound on its trace.

In this estimate, we condition on the realized stochastic Gram approximation and idealize the Monte Carlo gradient estimator as conditionally unbiased. Finite-chain MCMC bias and the plug-in error associated with the estimated normalization constant are not included explicitly; they may be represented as additional perturbation terms.

Let $\mathbf{d}_k$ denote the positive natural-gradient vector, the parameter update is $\boldsymbol{\theta}^{k+1}=\boldsymbol{\theta}^k-\tau_k\mathbf{d}_k$. If the regularized Gram system were solved exactly, then $\mathbf{H}_k\mathbf{d}_k=\widehat{\mathbf{g}}_k$. In the implemented matrix-free PCG solve, the residual $\boldsymbol{\rho}_k$ satisfies
\[
\mathbf{H}_k\mathbf{d}_k=\widehat{\mathbf{g}}_k-\boldsymbol{\rho}_k,
\qquad
\mathbf{d}_k=\mathbf{H}_k^{-1}\widehat{\mathbf{g}}_k-\mathbf{H}_k^{-1}\boldsymbol{\rho}_k .
\]
Because $\mathbf{G}_k$ is positive semidefinite, $\mathbf{H}_k\succeq \lambda_k\mathbf{I}$, and therefore
\[
\|\mathbf{H}_k^{-1}\|\le\frac1{\lambda_k},
\qquad
\|\mathbf{H}_k^{-1}\boldsymbol{\rho}_k\|\le\frac{\|\boldsymbol{\rho}_k\|}{\lambda_k}.
\]
Thus the same damping floor limits both Monte Carlo noise amplification and PCG residual amplification.

For the exact-solve idealization, set $\boldsymbol{\rho}_k=0$ and assume that $L$ is locally $L_{\mathrm{sm}}$-smooth in parameter space. The descent lemma, conditioned on $\mathcal{F}_k$, gives
\begin{equation}
\label{eq:sm-conditional-descent}
\begin{aligned}
\mathbb{E}[L(\boldsymbol{\theta}^{k+1})\mid\mathcal{F}_k]
&\le
L(\boldsymbol{\theta}^k)
-\tau_k \mathbf{g}_k^\top\mathbf{H}_k^{-1}\mathbf{g}_k  \\
&\qquad
+\frac{L_{\mathrm{sm}}\tau_k^2}{2}
\left(
\|\mathbf{H}_k^{-1}\mathbf{g}_k\|^2
+
\operatorname{tr}(\mathbf{H}_k^{-2}\boldsymbol{\Sigma}_k)
\right).
\end{aligned}
\end{equation}
If, additionally, $0\preceq \mathbf{G}_k\preceq M\mathbf{I}$ for a local spectral upper bound $M$, then
\begin{equation}
\label{eq:sm-spectral-descent}
\mathbb{E}[L(\boldsymbol{\theta}^{k+1})\mid\mathcal{F}_k]
\le
L(\boldsymbol{\theta}^k)
-\frac{\tau_k}{M+\lambda_k}\|\mathbf{g}_k\|^2
+\frac{L_{\mathrm{sm}}\tau_k^2}{2\lambda_k^2}
\left(\|\mathbf{g}_k\|^2+\sigma_k^2\right).
\end{equation}
The negative term is the expected descent contribution, whereas the positive term accounts for the finite step size and Monte Carlo noise. Increasing $\lambda_k$ suppresses the noise-amplification factor $1/\lambda_k^2$. The damping floor also prevents late-stage inverse-Gram amplification when $\mathbf{G}_k$ is singular or nearly singular.
For example, the conservative condition
\begin{equation}
\label{eq:sm-sufficient-descent}
\tau_k
\le
\frac{\lambda_k^2}{L_{\mathrm{sm}}(M+\lambda_k)}
\frac{\|\mathbf{g}_k\|^2}{\|\mathbf{g}_k\|^2+\sigma_k^2}
\end{equation}
is sufficient for expected descent under the exact-solve idealization above.

In the main experiments, we use the residual-aware schedule
\[
\lambda_k
=
\operatorname{clip}\!\left(10^{-3}\bar R_k,\,10^{-4},\,10^{-3}\right),
\]
where $\bar R_k$ is the exponentially averaged local variance. This quantity is a practical proxy for the noise and residual level, rather than an identity with $\sigma_k^2$. The resulting trajectory is shown in \Cref{subfig:lambda_lattice}. The schedule keeps stronger damping during the noisy early stage and retains a positive damping floor near convergence.

\section{Random Fourier feature ablations}
\label[appendix]{sm:rff-ablation}

\subsection{Fixed versus trainable Fourier matrix}
\label{sm:rff-fixed-trainable}
\leavevmode\par

To decide whether the Fourier matrix $\mathbf{B}$ should be fixed or trainable, we ran a controlled ablation on the 2D optical lattice benchmark with $\beta=250$, $m=32$, $\sigma=1.0$, and three random seeds.

\begin{table}[!htbp]
  \centering
  \caption{Fixed versus trainable Fourier matrix $\mathbf{B}$ for the 2D optical lattice benchmark ($m=32$, $\sigma=1.0$). Values are medians across three seeds; brackets denote interquartile ranges where reported.}
  \label{tab:supp-rff-trainable}
  \begin{tabular}{lccc}
    \toprule
    $\mathbf{B}$ & Best Res & Best $\varepsilon_{E}$ & $T_{\mathrm{epoch}}$ (s) \\
    \midrule
    Fixed     & 0.263 [0.172, 0.511] & $4.27\ [2.75, 5.20]\times10^{-5}$ & 0.155 \\
    Trainable & 0.083 [0.053, 0.099] & $4.63\ [3.34, 4.92]\times10^{-5}$ & 0.173 \\
    \bottomrule
  \end{tabular}
\end{table}

The energy error remains at the same order for both choices. Making $\mathbf{B}$ trainable reduces the best residual by roughly a factor of three, with only a modest increase in per-epoch cost. We therefore use trainable $\mathbf{B}$ in the optical-lattice experiments.

\subsection{Feature count}
\label{sm:rff-feature-count}
\leavevmode\par

We also swept the number of RFF features $m\in\{8,16,32,64\}$ with trainable $\mathbf{B}$ and $\sigma=1.0$ on the same 2D optical lattice benchmark.

\begin{table}[!htbp]
  \centering
  \caption{RFF feature-count ablation for the 2D optical lattice benchmark with trainable $\mathbf{B}$ and $\sigma=1.0$. Values are medians across three seeds; brackets denote interquartile ranges where reported.}
  \label{tab:supp-rff-count}
  \begin{tabular}{lccc}
    \toprule
    $m$ & Best Res & Best $\varepsilon_{E}$ & $T_{\mathrm{epoch}}$ (s) \\
    \midrule
    8  & 0.0208 [0.0121, 0.0213] & $2.93\ [2.21, 5.07]\times10^{-5}$ & 0.173 \\
    16 & 0.0231 [0.0175, 0.182]  & $1.22\ [1.10, 1.29]\times10^{-5}$ & 0.164 \\
    32 & 0.0831 [0.0531, 0.0994] & $4.63\ [3.34, 4.92]\times10^{-5}$ & 0.173 \\
    64 & 0.0416 [0.0282, 0.0792] & $7.12\ [4.89, 9.22]\times10^{-5}$ & 0.170 \\
    \bottomrule
  \end{tabular}
\end{table}

\Cref{tab:supp-rff-count} shows that the energy error is not sensitive to the feature count over this range. The choices $m=8$ and $m=16$ both give low residuals. Among the tested settings, $m=16$ gives the smallest median energy error and the lowest median epoch time. We therefore use $m=16$ for the optical-lattice experiments.

\subsection{High-dimensional RFF ablation}
\label{sm:rff-high-dimensional}
\leavevmode\par

The high-dimensional experiments raise a separate question: whether a fixed-size Gaussian RFF dictionary remains useful as the dimension increases. We tested this on the smooth $d=8$ anisotropic quadratic problem. The comparison includes raw coordinates without RFF, RFF with $m=16$, and RFF with $m=32$, using three seeds for each configuration.

\begin{figure}[!htbp]
  \centering
  \includegraphics[width=0.32\textwidth]{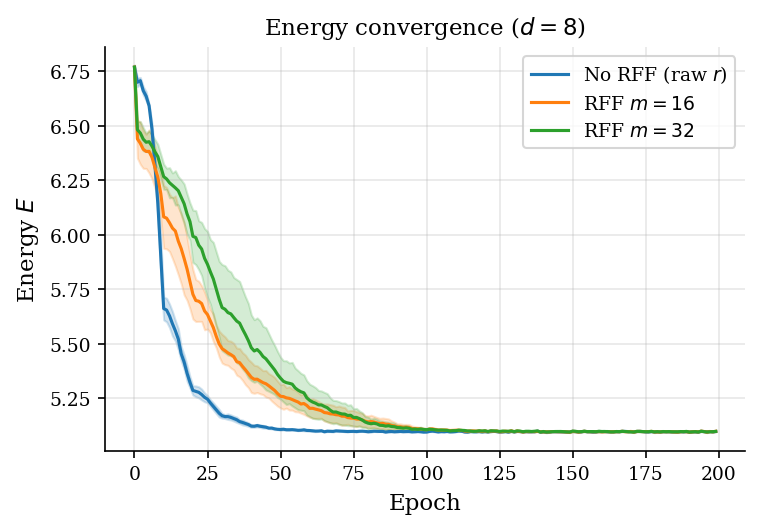}\hfill
  \includegraphics[width=0.32\textwidth]{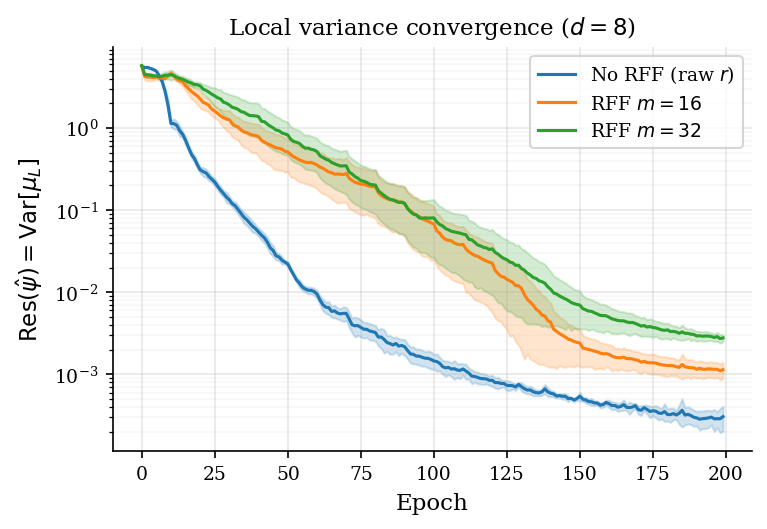}\hfill
  \includegraphics[width=0.32\textwidth]{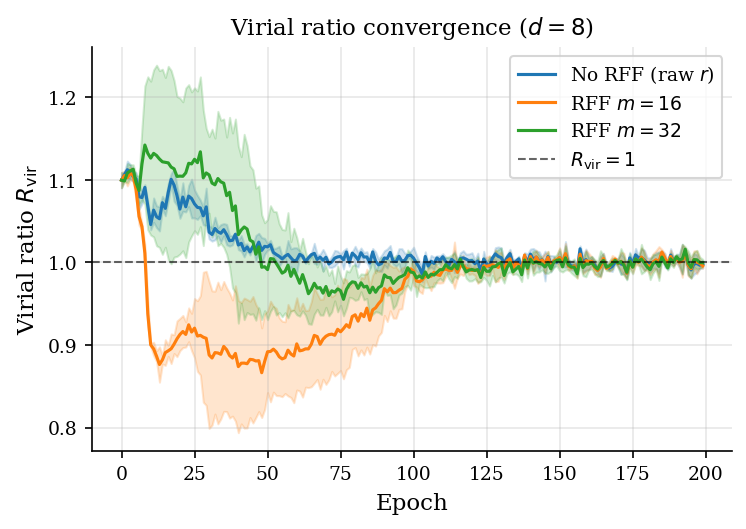}
  \caption{RFF ablation for the $d=8$ anisotropic quadratic potential. Left: energy convergence; middle: local variance on a logarithmic scale; right: virial ratio. Shaded bands denote one standard deviation across three seeds.}
  \label{fig:supp-rff-highdim}
\end{figure}

\begin{table}[!htbp]
  \centering
  \caption{Final statistics for the $d=8$ RFF ablation. Values are means $\pm$ standard deviations over three seeds.}
  \label{tab:supp-rff-highdim}
  \begin{tabular}{lccc}
    \toprule
    Configuration & Final $E$ & Final $\mathrm{Res}(\hat{\psi})$ & Final $R_{\mathrm{vir}}$ \\
    \midrule
    No RFF       & $5.098 \pm 0.001$ & $(3.1 \pm 1.0)\times10^{-4}$ & $0.998 \pm 0.004$ \\
    RFF ($m=16$) & $5.099 \pm 0.001$ & $(1.1 \pm 0.3)\times10^{-3}$ & $0.996 \pm 0.003$ \\
    RFF ($m=32$) & $5.098 \pm 0.001$ & $(2.8 \pm 0.3)\times10^{-3}$ & $1.000 \pm 0.004$ \\
    \bottomrule
  \end{tabular}
\end{table}

\Cref{fig:supp-rff-highdim,tab:supp-rff-highdim} show that the final energy and virial ratio are comparable across the three configurations. The local variance is lowest without RFF: $(3.1\pm1.0)\times10^{-4}$, compared with $(1.1\pm0.3)\times10^{-3}$ for $m=16$ and $(2.8\pm0.3)\times10^{-3}$ for $m=32$. Thus, for the smooth high-dimensional potential studied here, adding random Fourier features does not improve the result; instead, it increases the final residual. This behavior is consistent with the problem structure. The anisotropic quadratic potential produces a smooth, low-frequency ground state, so a high-frequency random dictionary is not the limiting factor and may complicate late-stage optimization.

\section{Nystr\"om preconditioning ablation}
\label{sm:nystrom-preconditioning-ablation}

\subsection{Nystr\"om rank}
\label[appendix]{sm:nystrom-rank}
\leavevmode\par

For the Nystr\"om rank study, we used the 2D optical lattice problem with $\beta=250$ and swept $r\in\{32,64,128,256,512\}$ over three random seeds. The Nystr\"om factor was refreshed every $K=10$ optimization steps. This ablation is reported for accuracy and spectral diagnostics. Timing results are omitted because a median computed only over non-refresh steps is not an amortized measure of the total optimization cost.

In this ablation, Hutch.\ rel.\ error denotes the mean probe-wise relative approximation error $\|Gz-MM^\top z\|_2/\|Gz\|_2$. Here $G$ is the matrix-free Gram matrix, and $MM^\top$ is the rank-$r$ Nystr\"om approximation. The quantity $\hat\kappa$ denotes the Lanczos estimate of the condition number of the damped operator after preconditioning, computed from a short PCG probe immediately after a Nystr\"om refresh.

\begin{table}[!htbp]
  \centering
  \caption{Nystr\"om rank ablation on the 2D optical lattice problem with $\beta=250$ and update frequency $K=10$. Values are medians across three seeds; brackets denote interquartile ranges.}
  \label{tab:supp-rank}
  \resizebox{\textwidth}{!}{
  \begin{tabular}{lccc}
    \toprule
    $r$ & Final $\mathrm{Res}(\hat\psi)$ & Hutch.\ rel. error & $\hat\kappa$ \\
    \midrule
    32 & 0.467 [0.269, 0.554] & 0.257 [0.224, 0.283] & $2.89{\times}10^4$ [$2.42{\times}10^4$, $3.29{\times}10^4$] \\
    64 & 0.195 [0.103, 0.213] & 0.184 [0.154, 0.204] & $1.02{\times}10^4$ [$7.72{\times}10^3$, $1.10{\times}10^4$] \\
    128 & 0.079 [0.050, 0.093] & 0.201 [0.154, 0.205] & $1.85{\times}10^3$ [$1.84{\times}10^3$, $2.71{\times}10^3$] \\
    256 & 0.078 [0.048, 0.208] & 0.181 [0.143, 0.193] & $1.67{\times}10^3$ [$1.61{\times}10^3$, $1.93{\times}10^3$] \\
    512 & 0.068 [0.041, 0.156] & 0.152 [0.132, 0.164] & $1.87{\times}10^3$ [$1.73{\times}10^3$, $2.37{\times}10^3$] \\
    \bottomrule
  \end{tabular}}
\end{table}

\begin{figure}[!htbp]
  \centering
  \includegraphics[width=0.47\linewidth]{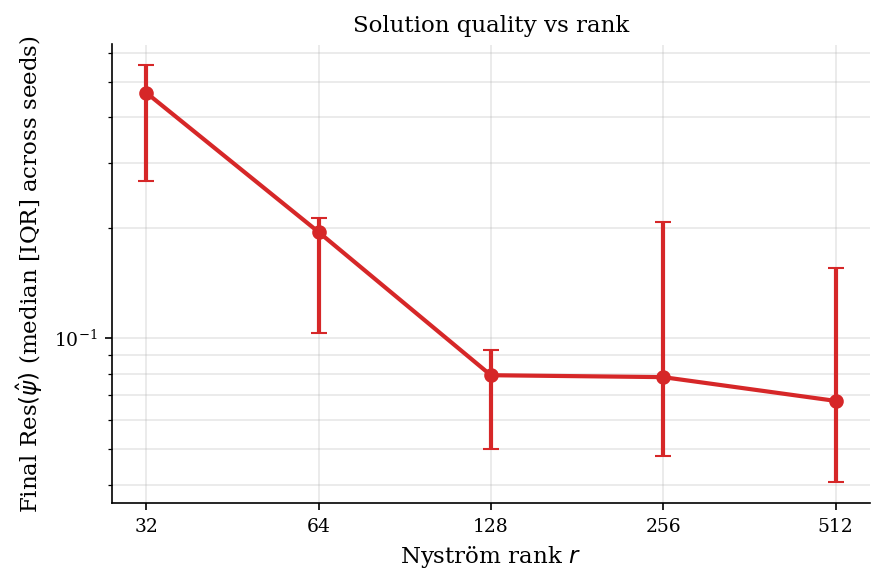}\hfill
  \includegraphics[width=0.47\linewidth]{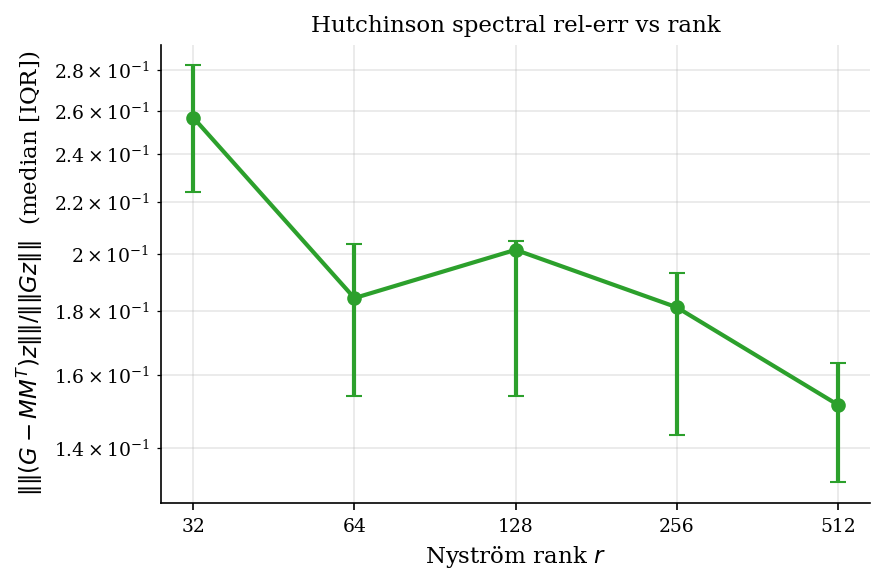}
  \caption{Nystr\"om rank ablation with update frequency $K=10$. Left: final residual; right: Hutchinson relative error.}
  \label{fig:supp-rank-sweep}
\end{figure}

\Cref{tab:supp-rank,fig:supp-rank-sweep} show a clear quality--rank tradeoff. Increasing $r$ from 32 to 128 reduces the median residual from 0.467 to 0.079 and lowers the Lanczos condition proxy from $2.89\times10^4$ to $1.85\times10^3$. The best median residual and Hutchinson relative error occur at $r=512$ (0.068 and 0.152, respectively), but the residual gain over $r=128$ is about 15\% while Nystr\"om basis storage grows linearly with $r$. We therefore use $r=128$ as a budget-aware default. The value $r=512$ is a quality-first option for comparable networks, whereas $r=32$ and $r=64$ are not recommended when high final accuracy is required.

\section{Adaptive defensive importance sampling}
\label{sm:adis}

\subsection{ADIS mixture parameter}
\label[appendix]{sm:adis-alpha}
\leavevmode\par

For the ADIS mixture parameter, we swept $\alpha\in\{0,0.2,0.5,0.8,0.95,1.0\}$ on two benchmarks: the 3D optical lattice benchmark, where a deterministic grid reference for $Z$ is available, and the $d=8$ anisotropic quadratic benchmark. Each configuration uses three independent seeds. We report both solution quality and importance-sampling diagnostics. Let $w_i=|\psi_{\boldsymbol\theta}(\mathbf r_i)|^2/q_\alpha(\mathbf r_i)$ be the raw ADIS importance weights, where $q_\alpha$ is the defensive proposal. The diagnostics use $N_{\mathrm{IS}}=8000$ samples in the 3D benchmark and $N_{\mathrm{IS}}=12{,}000$ in the $d=8$ benchmark. The ESS fraction is $\mathrm{ESS}/N_{\mathrm{IS}}$, where $\mathrm{ESS}=(\sum_i w_i)^2/\sum_i w_i^2$; larger values indicate less weight degeneracy. We use Pareto $\hat{k}$ only as a heuristic upper-tail diagnostic. Specifically, we sort the raw weights, retain the largest 20\%, set the threshold $u$ to the smallest retained weight, and fit a generalized Pareto distribution by maximum likelihood to the excesses $w_i-u$, with the location parameter fixed at zero; $\hat{k}$ is the fitted shape parameter. No Pareto smoothing is applied to the ADIS weights. Motivated by the Pareto-smoothed importance sampling \cite{vehtari2024pareto}, we regard $\hat{k}\geq0.7$ as evidence of a problematic upper tail at these sample sizes. A value below this empirical threshold is not, by itself, a guarantee of accuracy for the unsmoothed ADIS normalization estimate.

\begin{table}[!htbp]
  \centering
  \caption{ADIS $\alpha$ ablation for the 3D optical lattice problem with $\beta=200$. Values are medians across three seeds unless otherwise indicated.}
  \label{tab:supp-alpha-d3}
  \begin{tabular}{lcccc}
    \toprule
    $\alpha$ & ESS frac & Pareto $\hat{k}$ & $Z$ rel-err & Final Res \\
    \midrule
    0.0  & 0.032 & \textbf{3.82} & 0.055 & 0.029 \\
    0.2  & 0.182 & $-1.25$       & 0.0053 & 0.0041 \\
    0.5  & 0.385 & $-0.56$       & 0.0033 & 0.0021 \\
    \textbf{0.8} & 0.589 & $-0.16$ & $\mathbf{0.0013}$ & 0.0024 \\
    0.95 & 0.688 & $0.02$        & 0.0029 & 0.0028 \\
    1.0  & 0.720 & $0.01$        & 0.0075 & 0.0035 \\
    \bottomrule
  \end{tabular}
\end{table}

\begin{table}[!htbp]
  \centering
  \caption{ADIS $\alpha$ ablation for the $d=8$ anisotropic quadratic problem with $\beta=2000$. Values are medians across three seeds unless otherwise indicated.}
  \label{tab:supp-alpha-d8}
  \begin{tabular}{lccc}
    \toprule
    $\alpha$ & ESS frac & Pareto $\hat{k}$ & Final Res \\
    \midrule
    0.0  & 0.0001 & \textbf{8.93}  & 0.083 \\
    0.2  & 0.099  & $0.02$         & 0.00049 \\
    0.5  & 0.249  & $-0.37$        & 0.00046 \\
    \textbf{0.8} & 0.394 & $-0.39$ & 0.00045 \\
    0.95 & 0.467  & $-0.39$        & \textbf{0.00044} \\
    1.0  & 0.497  & $-0.41$        & 0.00048 \\
    \bottomrule
  \end{tabular}
\end{table}

\Cref{tab:supp-alpha-d3,tab:supp-alpha-d8} reveal complementary failure modes at the two endpoints. The pure uniform proposal, $\alpha=0$, has a severely problematic upper-weight tail, with $\hat{k}=3.82$ in $d=3$ and $\hat{k}=8.93$ in $d=8$. In $d=8$, its final residual is 0.083, about $185\times$ larger than that for $\alpha=0.8$. The pure Gaussian proposal, $\alpha=1$, has a high ESS fraction but gives a worse $Z$ estimate in the $d=3$ benchmark. Its $Z$ relative error is $7.5\times10^{-3}$, compared with $1.3\times10^{-3}$ for $\alpha=0.8$. Thus ESS alone is not sufficient; the tail diagnostic, reference normalization error when available, and final residual should be considered together. Both components of the defensive mixture therefore play a useful role. We use $\alpha=0.8$ as the conservative default because it gives the best $Z$ relative error in the $d=3$ benchmark, has below-threshold Pareto diagnostics in both dimensions ($\hat{k}=-0.16$ and $-0.39$), and is within about 2\% of the best median final residual in the $d=8$ sweep. For new problems, we recommend monitoring ESS and Pareto $\hat{k}$ and, when possible, checking values in the range $\alpha\in[0.5,0.95]$.

\bibliographystyle{siamplain}
\bibliography{refs}
\end{document}